\documentclass[12pt,oneside,reqno]{amsart}

\usepackage{txfonts}
\usepackage{bbm}
\usepackage{amsmath}
\usepackage{graphicx}
\usepackage{mathrsfs}
\usepackage{stmaryrd}
\usepackage{amsfonts}
\usepackage{enumerate,amsmath,amssymb,amsthm}

\numberwithin{equation}{section}

\newcommand{\be}{\begin{eqnarray}}
\newcommand{\ee}{\end{eqnarray}}
\newcommand{\ce}{\begin{eqnarray*}}
\newcommand{\de}{\end{eqnarray*}}
\newtheorem{theorem}{Theorem}[section]
\newtheorem{lemma}[theorem]{Lemma}
\newtheorem{remark}[theorem]{Remark}
\newtheorem{definition}[theorem]{Definition}
\newtheorem{proposition}[theorem]{Proposition}
\newtheorem{Examples}[theorem]{Example}
\newtheorem{corollary}[theorem]{Corollary}

\def\eps{\varepsilon}

\def\p{\partial}

\def\[{{\Big[}}
\def\]{{\Big]}}
\def\<{{\langle}}
\def\>{{\rangle}}
\def\({{\Big(}}
\def\){{\Big)}}

\def\bx{{\mathbf{x}}}

\def\dif{{\mathord{{\rm d}}}}

\def\no{\nonumber}
\def\={&\!\!=\!\!&}
\def\bt{\begin{theorem}}
\def\et{\end{theorem}}
\def\bl{\begin{lemma}}
\def\el{\end{lemma}}
\def\br{\begin{remark}}
\def\er{\end{remark}}
\def\bd{\begin{definition}}
\def\ed{\end{definition}}
\def\bp{\begin{proposition}}
\def\ep{\end{proposition}}
\def\bc{\begin{corollary}}
\def\ec{\end{corollary}}
\def\bx{\begin{Examples}}
\def\ex{\end{Examples}}

\def\cE{{\mathcal E}}

\def\cJ{{\mathcal J}}

\def\cL{{\mathcal L}}

\def\cP{{\mathcal P}}

\def\cS{{\mathcal S}}

\def\mD{{\mathbb D}}
\def\mE{{\mathbb E}}

\def\mH{{\mathbb H}}
\def\mI{{\mathbb I}}

\def\mN{{\mathbb N}}

\def\mP{{\mathbb P}}

\def\mR{{\mathbb R}}
\def\mS{{\mathbb S}}

\def\mU{{\mathbb U}}

\def\mW{{\mathbb W}}

\def\sA{{\mathscr A}}
\def\sB{{\mathscr B}}
\def\sC{{\mathscr C}}

\def\sF{{\mathscr F}}

\def\sI{{\mathscr I}}

\def\sS{{\mathscr S}}

\def\geq{\geqslant}
\def\leq{\leqslant}

\def\div{\mathord{{\rm div}}}

\def\v{{\mathrm v}}
\def\e{{\mathrm e}}

\allowdisplaybreaks

\begin{document}

\title{Nonlocal H\"ormander's hypoellipticity theorem}

\date{}
\author{XICHENG ZHANG
\\
\\
{\it Dedicated to the memory of Professor Paul Malliavin}}

%\thanks{{\it Keywords: }Reisz transform, critical parabolic equation}

%\thanks{$*$This work is supported by NSFs of China (No. 10971076) and Program for New Century Excellent Talents in University (NCET 10-0654).}

\address{
School of Mathematics and Statistics,
Wuhan University, Wuhan, Hubei 430072, P.R.China\\
Email: XichengZhang@gmail.com
 }

\begin{abstract}
Consider the following nonlocal integro-differential operator: for $\alpha\in(0,2)$,
$$
\cL^{(\alpha)}_{\sigma,b} f(x):=\mbox{p.v.} \int_{\mR^d-\{0\}}\frac{f(x+\sigma(x)z)-f(x)}{|z|^{d+\alpha}}\dif z+b(x)\cdot\nabla f(x),
$$
where $\sigma:\mR^d\to\mR^d\times\mR^d$ and $b:\mR^d\to\mR^d$ are two $C^\infty_b$-functions, and p.v. stands for the Cauchy principal value.  
Let $B_1(x):=\sigma(x)$ and $B_{j+1}(x):=b(x)\cdot\nabla B_j(x)-\nabla b(x)\cdot B_j(x)$ for $j\in\mN$. Under the following H\"ormander's type condition:
for any $x\in\mR^d$ and some $n=n(x)\in\mN$,
$$
\mathrm{Rank}[B_1(x), B_2(x),\cdots, B_n(x)]=d,
$$
by using the Malliavin calculus, we prove the existence of the heat kernel $\rho_t(x,y)$ to the operator $\cL^{(\alpha)}_{\sigma,b}$ 
as well as the continuity of $x\mapsto \rho_t(x,\cdot)$ in $L^1(\mR^d)$ for each $t>0$.
Moreover, when $\sigma(x)=\sigma$ is constant, under the following uniform H\"ormander's type condition: for some $j_0\in\mN$,
$$
\inf_{x\in\mR^d}\inf_{|u|=1}\sum_{j=1}^{j_0}|u B_j(x)|^2>0,
$$
we also prove the smoothness of $(t,x,y)\mapsto\rho_t(x,y)$ with $\rho_t(\cdot,\cdot)\in C^\infty_b(\mR^d\times\mR^d)$ for each $t>0$.
\end{abstract}

\maketitle \rm

\section{Introduction}
\subsection{Background}
Let $A_0, A_1,\cdots, A_d$ be $d+1$-smooth vector fields on $\mR^d$. Consider the following second order differential operator
\begin{align}
\cL^{(2)}:=\sum_{k=1}^d A_k^2+A_0.\label{Se}
\end{align}
Let $\sA_0:=\{A_1,\cdots,A_d\}$ and for $j\in\mN$,
\begin{align}
\sA_j:=\Big\{[A, A_k]: A\in\sA_{j-1}, k=0,\cdots,d\Big\},\label{Lie}
\end{align}
where $[A, A_i]:=AA_i-A_iA$ is the usual Lie bracket. The celebrated H\"ormander's hypoellipticity theorem asserts that if 
\begin{align}
\mbox{$\cup_{j=0}^\infty\sA_j$ spans $\mR^d$ at each point $x$,}\label{Hor}
\end{align}
then $\p_t-\cL^{(2)}$ has a smooth fundamental solution. In 1967,
H\"ormander \cite{Ho} first proved this result by using the analytic argument. Later, in 1976, Malliavin \cite{Ma} provided a probabilistic proof. 
Malliavin's idea was further developed by Bismut \cite{Bi},  Kusuoka-Stroock \cite{Ku-St1, Ku-St2}, Watanabe \cite{Wa}, Norris \cite{No} etc.
Notice that H\"ormander's condition (\ref{Hor}) is a local condition in the sense that if $\cup_{j=0}^{n(x)}\sA_j$ spans $\mR^d$ at point $x$, then
$\cup_{j=0}^{n(x)}\sA_j$ also spans $\mR^d$ at a neighbourhood of $x$.

In contrary to the {\it local} operator (\ref{Se}),  we consider the following {\it nonlocal} operator
\begin{align}
\cL^{(\alpha)}_{\sigma,b} f(x):=\mbox{p.v.} \int_{\mR^d-\{0\}}\frac{f(x+\sigma(x)z)-f(x)}{|z|^{d+\alpha}}\dif z+b(x)\cdot\nabla f(x),\ \ \alpha\in(0,2),\label{OP}
\end{align}
where $\sigma:\mR^d\to\mR^d\times\mR^d$ and $b:\mR^d\to\mR^d$ are two smooth functions, and p.v. stands for the Cauchy principal value. Notice that 
if $\sigma(x)=\mI$ is the identity matrix and $b=0$, then
$$
\cL^{(\alpha)}_{\mI,0} f=c_{d,\alpha}\Delta^{\frac{\alpha}{2}}f,
$$
where $\Delta^{\frac{\alpha}{2}}$ is the usual fractional Laplacian and $c_{d,\alpha}$ is a constant. 
Recently, there are increasing interests of studying such type of integro-differential operators. For example, in  \cite{Mo-Xu},  
Morimoto and Xu first studied the hypoelliptic effects to the following fractional kinetic Fokker-Planck equation
\begin{align}
\p_t f(t,x,\v)+\v\cdot\nabla_x f(t,x,\v)=a(x,\v)\Delta^{\frac{\alpha}{2}}_\v f(t,x,\v)+g,\label{OP2}
\end{align}
where $a$ is a bounded smooth and strictly positive function. See also \cite{Al-Mo-Uk-Xu-Ya, Ch-Li-Xu, Al} etc. for some related works about this equation. 
The physical motivation of studying equation (\ref{OP2}) is refereed to \cite{Za}.
Another motivation is linked with a linear model of spatially inhomogeneous Boltzmann equation without an angular cutoff (see \cite{Mo-Xu, Al-Mo-Uk-Xu-Ya}). 
 Notice that the operator in equation (\ref{OP2}) is degenerate, and can be written as (\ref{OP}) with nonlinear $\sigma$ and linear $b$ by the change of variables.

On the other hand, let us consider the following stochastic nonlinear oscillators studied in \cite{Ec-Ha, Re-Th, Ca0} etc.:
\begin{align}\label{6-1}
\begin{cases}
\dif z_i(t)=u_i(t)\dif t,& i=1,\cdots,d,\\
\dif u_i(t)=-\partial_{z_i}H(z(t),u(t))\dif t, & i=2,\cdots,d-1,\\ 
\dif u_i(t)=-[\partial_{z_i}H(z(t),u(t))+\gamma_iu_i(t)]\dif t+\sqrt{T_i}\dif L^i_t, & i=1, d,
\end{cases}
\end{align}
where $L^1_t$ and $L^d_t$ are two independent of $\alpha$-stable processes, $d\geq 3$, $\gamma_1,\gamma_d\in\mR$, $T_1,T_d>0$, and
$$
H(z,u):=\sum_{i=1}^d\left(\frac{1}{2}|u_i|^2+V(z_i)\right)+\sum_{i=1}^{d-1}U(z_{i+1}-z_{i}).
$$
The typical examples of potential functions $V$ and $U$ are
$$
V(z)=\frac{|z|^2}{2},\ \ U(z)=\frac{|z|^2}{2}+\frac{|z|^4}{4}.
$$
The Hamiltonian $H$ describes a chain of particles with nearest-neighbor interaction.  We remark that the generator of stochastic equation (\ref{6-1})
also takes the form (\ref{OP}) with constant $\sigma$ and nonlinear $b$, which is a typical H\"ormander's type operator if we assume $U''>0$ (cf. \cite{Do-Pe-So-Zh}).

In the present work, we are interested in the following problem: Under what degenerate conditions on $\sigma$ and $b$ appearing in (\ref{OP}), 
we can prove a similar theorem for integro-differential operator $\cL^{(\alpha)}_{\sigma,b}$ as in H\"ormander's theorem.
The purpose of this paper is to give partial results under H\"ormander's type conditions.
\subsection{Main results}
Let $B_1(x):=\sigma(x)$ and  for $j\in\mN$, define
\begin{align}
B_{j+1}(x):=b(x)\cdot\nabla B_j(x)-B_j(x)\cdot\nabla b(x).\label{EUT1}
\end{align}
Our first result is the following existence of the heat kernel to $\cL^{(\alpha)}_{\sigma,b}$ and the $L^1$-continuity of the heat kernel with respect to the first variable.
\bt\label{Th00}
Let $b:\mR^d\rightarrow\mR^d$ and $\sigma:\mR^d\rightarrow\mR^d\times\mR^d$ be $C^{\infty}$-functions 
with bounded partial derivatives of first order. Assume that for each $x\in\mR^d$ and some $n=n(x)\in\mN$,
\begin{align}
\mathrm{Rank}[B_1(x), B_2(x),\cdots, B_n(x)]=d,\label{ForH}
\end{align}
then there exists a function $\rho_t(x,y)$ such that for any $f\in C^2_b(\mR^d)$ and $t>0, x\in\mR^d$,
\begin{align}
\p_t\int_{\mR^d}f(y)\rho_t(x,y)\dif y=\int_{\mR^d}\cL^{(\alpha)}_{\sigma,b}f(y)\rho_t(x,y)\dif y.\label{RT1}
\end{align}
Moreover, $(t,x)\mapsto\rho_t(x,\cdot)$ is continuous in $L^1(\mR^d)$ on $(0,\infty)\times\mR^d$.
In particular, the semigroup defined by $\cP_t f(x):=\int_{\mR^d}f(y)\rho_t(x,y)\dif y$ has the strong Feller property.
\et
\br
If $\sigma(x)=\sigma$ is constant, in \cite{Zh2} and \cite{Do-Pe-So-Zh}, we have already proved this theorem.
Moreover, if we define vector fields
$$
A_0(x):=\sum_{i=1}^db_i(x)\p_i,\ \ A_k(x):=\sum_{i=1}^d\sigma_{ki}(x)\p_i,
$$
and set $\sA_0:=\{A_1,\cdots, A_d\}$ and
$$
\sA_j:=\Big\{[A,A_0]: A\in\sA_{j-1}\Big\},\ j\in\mN,
$$
then (\ref{ForH}) is equivalent that
$$
\cup_{j=0}^\infty\sA_j\mbox{ spans $\mR^d$ at each point $x\in\mR^d$}.
$$
Compared with (\ref{Hor}), in the definition of $\sA_j$, all the Lie bracket is taken only with the drift vector field $A_0$.
\er

About the smoothness of $\rho_t(x,y)$, we have the following partial result.
\bt\label{Th1}
Assume that $\sigma(x)=\sigma$ is constant, $b$ has bounded partial derivatives of all orders and for some $j_0\in\mN$,
\begin{align}
\inf_{x\in\mR^d}\inf_{|u|=1}\sum_{j=1}^{j_0}|u B_j(x)|^2=:\kappa_1>0.\label{UH}
\end{align}
Then there exists a probability density function $\rho_t(x,y)\in C^\infty(\mR_+\times\mR^d\times\mR^d)$ with
\begin{align}
\rho_t(\cdot,\cdot)\in C^\infty_b(\mR^d\times\mR^d)\ \forall t>0,\label{ET9}
\end{align}
and such that
\begin{align}
\p_t\rho_t(x,y)=\cL^{(\alpha)}_{\sigma,b}\ \rho_t(\cdot,y)(x)\ \forall t>0,\label{ET8}
\end{align}
with $\lim_{t\downarrow 0}\rho(t,x,y)=\delta_y(x)$ in the distributional sense.
\et
\br
In the linear case of $b(x)=Bx$, Priola and Zabczyk \cite{Pr-Za} first established this result under (\ref{UH}), and Kulik \cite{Ku00} also showed that
(\ref{UH}) is necessary for the existence of smooth fundamental solutions (see also \cite{Si} for the related studies). In the nonlinear case, when $j_0=2$, this
result was proven in \cite{Zh2}. Here an open question is the smoothness of $\rho_t(x,y)$ for the general non-constant $\sigma(x)$. 
We will study this problem in a subsequent work \cite{Zh5}.
\er
\subsection{Arguments and Related results}
Our proofs of Theorems \ref{Th00} and \ref{Th1} are based on the Malliavin calculus. Let $(L_t)_{t\geq 0}$ 
be a rotationally invariant $\alpha$-stable process. Consider the following
stochastic differential equation (abbreviated as SDE) driven by $L_t$:
$$
\dif X_t(x)=b(X_t(x))\dif t+\sigma(X_{t-}(x))\dif L_t,\ \ X_0(x)=x.
$$
It is well known that the generator of this SDE is given by $\cL^{(\alpha)}_{\sigma,b}$. 
Our task is thus to show that $X_t(x)$ admits a $L^1$-continuity density $\rho_t(x,y)$ under (\ref{ForH}), and which is smooth under (\ref{UH}) when $\sigma(x)=\sigma$
is constant. Let $(W_t)_{t\geq 0}$ be a $d$-dimensional standard Brownian motion and $(S_t)_{t\geq 0}$ an independent $\alpha/2$-stable subordinator. 
Then $(W_{S_t})_{t\geq 0}$ (also called subordinated Brownian motion) is a rotationally invariant $\alpha$-stable process, and we shall assume
$$
L_t=W_{S_t}.
$$
As in \cite{Ku, Zh1, Zh2}, by taking regular conditional probability with respect to $S_\cdot$, the solution $X_t(x)$ can be regarded as a Brownian functional so that we can use the classical
Malliavin calculus to study the density of $X_t(x)$. 

Let us now recall some well known results about the smooth density for SDEs with jumps.
In \cite{Bi0}, Bismut first introduced the analogue of the Malliavin calculus for SDEs with jumps by the perturbation of jump sizes and using Girsanov's transformation.
In \cite{Bi-Ja-Gr}, Bichteler, Gravereaux and Jacod systematically studied  the smooth density for SDEs driven by nondegenerate jump noises.
In \cite{Pi}, Picard used the difference operator to give another criterion for the smoothness of 
the distribution density of Poisson functionals, and also applied it to SDEs driven by pure jump L\'evy processes. 
By combining the classical Malliavin calculus and Picard's difference operator argument, 
Ishikawa and Kunita in \cite{Is-Ku} obtained a new criterion for the smooth density of Wiener-Poisson functionals. In \cite{Ku}, Kusuoka developed the Malliavin calculus for subordinated Brownian functionals as mentioned above. The advantage of this method lies in the fact that it is not necessary to develop a new Malliavin calculus for jump processes
(see \cite{ Wa-Wa, Zh1, Zh2, Zh3} for more applications). 
It is noticed that in all these results, the jump noises of SDEs are assumed to be nondegenerate.
On the other hand, Cass  in \cite{Ca} established a H\"ormander's type theorem for SDEs with jumps
by proving a Norris' type lemma for discontinuous semimartingales,  but the Brownian diffusion term can not disappear. 
In the pure jump degenerate case, by using a Komatsu-Takeuchi's estimate proven in \cite{Ko-Ta1} for discontinuous semimartingales, 
Takeuchi \cite{Ta} and Kunita \cite{Ku01, Ku0} 
also obtained similar H\"ormander's theorems. Unfortunately, the estimate in \cite{Ko-Ta1} seems incomplete, and which will be discussed below.

As usual, in the proof of Theorem \ref{Th00},  we only need to check the invertibility of  the Malliavin covariance matrix $\Sigma_t(x)$
defined by (\ref{Ep4}) below. To obtain the $L^1$-continuity of the density, we shall use a criterion proved by 
Bogachev \cite[Corollary 9.6.12]{Bo} about the convergence of a family of Wiener functionals with respect to the total variation distance
(see also \cite{Ma-Po, Ba-Ca, Do-Pe-So-Zh}). While in the proof of Theorem \ref{Th1},
the main difficulty is to prove the $L^p$-integrability of the inverse of $\Sigma_t(x)$, 
where the key point is to prove a Norris'  type lemma for discontinuous semimartingales (cf. \cite{Ca}). 
In the continuous diffusion case, Komatsu-Takeuchi \cite{Ko-Ta0}
used a new estimate for continuous semimartingales to give a simplified probabilistic proof for H\"ormander's theorem.
Later, in \cite{Ko-Ta1}, they also extended it to the case of discontinuous semimartingales so that they can prove a H\"ormander's type theorem for
the SDEs driven by general jump processes. However, it seems that there is a gap in their proofs.
Here, we shall provide a slightly different and simplified proof for their estimate.  In particular, the parameters appearing in the estimate
are explicitly computed, which is crucial for our proof.
As in \cite{Ca}, the difference with \cite{Ko-Ta1} is that the jump size of the discontinuous part has to be small.

\subsection{Layout and Notations}
This paper is organised as follows: In Section 2, we give some preliminaries about the Malliavin calculus and some necessary lemmas for later use.
In particular, the Malliavin covariance matrix is calculated.
In Section 3, we prove Theorem \ref{Th00}. In Section 4, we prove a Komatsu-Takeuchi's estimate for discontinuous semimartingales 
by using different calculations. In Section 5, Theorem \ref{Th1} will be proven.

Before concluding this introduction, we collect some notations and make some conventions for later use.
\begin{itemize}
\item Write $\mR_+=(0,\infty)$ and $\mN_0:=\mN\cup\{0\}$, $\mR^d_0:=\mR^d-\{0\}$.
\item For a c\`adl\`ag function $f:\mR_+\to\mR^d$, $\Delta f_s:=f_s-f_{s-}$.
\item The inner product in Euclidean space is denoted by $\<x,y\>$ or $x\cdot y$.
\item For $p>0$, $\Gamma(p):=\int^\infty_0 s^{p-1} \e^{-s}\dif s$ denotes the usual Gamma function.
\item $\cS(\mR^d)$: The Schwardz space of rapidly decreasing  smooth functions.
\item $C^\infty_b(\mR^d)$: The space of all smooth bounded functions with bounded partial derivatives of all orders.
\item $\nabla=(\p_1,\cdots,\p_d)$ denotes the gradient operator, and $D$ the Malliavin derivative operator.
\item For a smooth function $f:\mR^d\to\mR^d$, $(\nabla f)_{ij}:=(\p_j f^i)$ denotes the Jacobian matrix of $f$.
\item The asterisk $*$ denotes the transpose of a matrix or a column vector, or the dual operator.
\item The capital letter $C$ (resp. the small letter $c$) with or without index will denote an unimportant constant with values in $[1,\infty)$ (resp. $(0,1)$).
\end{itemize}

\section{Preliminaries}
We first introduce the canonical space of a subordinated Brownian motion $W_{S_t}$. Let $(\mW,\mH,\mu_\mW)$ be the classical Wiener space, i.e., 
$\mW$ is the space of all continuous functions from $\mR_+$ to $\mR^d$ with vanishing values at starting point $0$, $\mH\subset\mW$ is the Cameron-Martin space consisting of all
absolutely continuous functions with square integrable derivatives, $\mu_\mW$ is the Wiener measure so that the coordinate process
$$
W_t(w):=w_t
$$
is a $d$-dimensional standard Brownian motion. Let $\mS$ be the space of all c\`adl\`ag increasing functions from $\mR_+$ to $\mR_+$ with
$\ell_0=0$. Suppose that $\mS$ is endowed with the Skorohod metric and a probability measure $\mu_\mS$ so that the coordinate process
$$
S_t(\ell):=\ell_t
$$
is a pure jump subordinator with L\'evy measure $\nu_S$ satisfying 
$$
\int^\infty_0(1\wedge u)\nu_S(\dif u)<+\infty.
$$
Consider now the following product probability space
$$
(\Omega,\sF,\mP):=\Big(\mW\times \mS, \sB(\mW)\times\sB(\mS), \mu_\mW\times\mu_{\mS}\Big),
$$
and define for $\omega=(w,\ell)\in\mW\times \mS$,
$$
L_t(\omega):=w_{\ell_t}.
$$
Then $(L_t)_{t\geq 0}$ is a $d$-dimensional pure jump L\'evy process with L\'evy measure $\nu_L$ given by
\begin{align}
\nu_L(E)=\int^\infty_0\!\!\int_E(2\pi u)^{-\frac{d}{2}}\e^{-\frac{|z|^2}{2u}}\dif z\nu_S(\dif u),\ \ E\in\sB(\mR^d).\label{Levy}
\end{align}
In particular, if $\nu_S(\dif u)=u^{-1-\frac{\alpha}{2}}\dif u$, then
$$
\nu_L(\dif z)=2^{-\frac{d+\alpha}{2}}(2\pi)^{-\frac{d}{2}}\Gamma(\tfrac{d+\alpha}{2})|z|^{-d-\alpha}\dif z.
$$

\subsection{Conditional Malliavin calculus on $\Omega$}
Let us now recall some basic notions and facts about the Malliavin calculus (cf. \cite{Nu}).
Let $\mU$ be a real separable Hilbert space. Let $\sC(\mW;\mU)$ be the class of all $\mU$-valued smooth cylindrical functionals on $\mW$ with the form:
$$
F=\sum_{i=1}^mf_i(W(h_1),\cdots,W(h_n)) u_i,
$$
where $f_i\in C^\infty_b(\mR^n)$, $u_i\in\mU$, $h_1,\cdots, h_n\in\mH$, and
$$
W(h):=\sum_{k=1}^d\int^\infty_0 \dot h^k_s\dif W^k_{s}.
$$
The Malliavin derivative of $F$ is defined by
$$
DF:=\sum_{i=1}^m\sum_{j=1}^n(\p_j f_i)(W(h_1),\cdots,W(h_n))u_i\otimes h_j\in\mU\otimes\mH.
$$
By an iteration argument, for any $k\in\mN$, the higher order Malliavin derivative $D^kF$ of $F$ can be defined as a random variable in $\mU\otimes\mH^{\otimes k}$. 
It is well known that the operator $(D^k,\sC(\mW;\mU))$ is closable from $L^p(\mW;\mU)$ to $L^p(\mW;\mU\otimes\mH^{\otimes k})$ for each $p\geq 1$ 
(cf. \cite[p.26, Proposition 1.2.1]{Nu}). For every  $p\geq 1$ and $k\in\mN$, we introduce a norm on $\sC(\mW;\mU)$ by
$$
\|F\|_{k,p}:=\left(\mE^{\mu_\mW}|F|^p+\sum_{l=1}^k\mE^{\mu_\mW}\(\|D^lF\|^p_{\mU\otimes\mH^{\otimes l}}\)\right)^{\frac{1}{p}}.
$$
The Wiener-Sobolev space $\mD^{k,p}(\mW;\mU)$ is defined as the closure of $\sC(\mW;\mU)$ with respect to the above norm.
The dual operator $D^*$ of $D$ (also called divergence operator) is defined by relation
$$
\mE^{\mu_\mW} \<DF,U\>_\mH=\mE^{\mu_\mW} (FD^*U),\ \ U\in \mD^{1,2}(\mW;\mH).
$$
The following Meyer's inequality holds (cf. \cite[p.75, Proposition 1.5.4]{Nu}). For any $p>1$ and $U\in \mD^{1,p}(\mW;\mH)$,
\begin{align}
\|D^* U\|_p\leq C_p \|U\|_{1,p}.\label{Me}
\end{align}

Now, let $F(w,\ell)$ be a $\mU$-valued functional on $\Omega=\mW\times\mS$. Let us denote 
$$
F^\ell(w):=F(w,\ell),
$$
and define
\begin{align*}
\mD^{k,p}(\mU):=\mD^{k,p}(\Omega;\mU)&:=\Bigg\{F\in L^p(\Omega;\mU): \mbox{for $\mu_\mS$-a.a. $\ell\in\mS$, }F^\ell\in \mD^{k,p}(\mW;\mU),\\
&\qquad\mbox{ and }\|F\|_{k,p}:=\int_{\Omega}\|D^k F^\ell(w)\|^p_{\mU\otimes\mH^{\otimes k}}\mP(\dif w\times\dif\ell)<\infty\Bigg\}
\end{align*}
and
\begin{align}
(D^kF)(w,\ell):=D^kF^\ell(w), \ \ F\in \mD^{k,p}(\mU).\label{Def1}
\end{align}
Below we simply write
$$
\mD^\infty(\mU):=\mD^\infty(\mU):=\cap_{k\in\mN,p\geq 1}\mD^{k,p}(\mU)
$$
and
$$
\mD^{k,p}:=\mD^{k,p}(\mR^d),\ \ \mD^\infty:=\mD^\infty(\mR^d).
$$
Let $F=(F_1,\cdots, F_d)$ be a random vector on $\Omega$ with $F^\ell\in\mD^{1,2}(\mW;\mR^d)$ for each $\ell\in\mS$. The Malliavin covariance matrix of $F$ is defined by
\begin{align}
(\Sigma_F)_{ij}(w,\ell):=\<DF^\ell_i(w), DF^\ell_j(w)\>_\mH.\label{Def2}
\end{align}

The following criterion about the $L^1$-continuity of the densities of a random field can be found  in \cite[Corollary 9.6.12]{Bo} (see also \cite{Ma-Po, Ba-Ca, Do-Pe-So-Zh}).
\bt\label{Th11}
Let $\Lambda$ be a metric space and $(X_\lambda)_{\lambda\in\Lambda}$ a family of  $\mR^d$-valued functionals over $\Omega$. Suppose that
for some $p\geq d$,
\begin{enumerate}[{\bf(H1)}]
\item For each $\lambda\in\Lambda$ and $\mu_\mS$-almost all $\ell\in\mS$, 
$X^\ell_\lambda\in \mD^{1,p}(\mW;\mR^d)$, and if $\lambda_n\to\lambda$, then
$$
\lim_{n\to\infty}\|X^\ell_{\lambda_n}-X^\ell_\lambda\|_{1,p}=0.
$$
\item For each $\lambda\in\Lambda$, the Malliavin covariance matrix $\Sigma_{X_\lambda}$ of $X_\lambda$ is invertible almost surely.
\end{enumerate}
Then the law of $X_\lambda$ in $\mR^d$ admits a density $\rho_\lambda(x)$ so that
$\lambda\mapsto\rho_\lambda(\cdot)$ is continuous in $L^1(\mR^d)$.
\et
\begin{proof}
First of all, by {\bf (H2)} and Fubini's theorem, for $\mu_\mS$-almost all $\ell\in\mS$, $\Sigma_{X^\ell_{\lambda}}(w)$ 
is invertible for $\mu_\mW$-almost all $w$. Thus, by {\bf (H1)} and Bouleau-Hirsch's criterion \cite[p.92, Theorem 2.1.1]{Nu},
the law of $X^\ell_\lambda$ is absolutely continuous with respect to the Lebesgue measure, and so does $X_\lambda$. The density is denoted by $\rho_\lambda(x)$.
Let $\lambda_n\to\lambda$. By {\bf (H1)} and Bogachev's criterion \cite[Corollary 9.6.12]{Bo} (see also \cite[Corollary 2.2]{Ma-Po}), 
we obtain that for $\mu_\mS$-almost all $\ell\in\mS$,
$$
\lim_{n\to\infty}\sup_{\|f\|_\infty\leq 1}|\mE^{\mu_\mW}(f(X^\ell_{\lambda_n})-f(X^\ell_{\lambda}))|=0,
$$
which, by the dominated convergence theorem, implies that 
\begin{align}
\lim_{n\to\infty}\sup_{\|f\|_\infty\leq 1}|\mE(f(X_{\lambda_n})-f(X_{\lambda}))|
\leq \int_{\mS}\lim_{n\to\infty}\sup_{\|f\|_\infty\leq 1}|\mE^{\mu_\mW}(f(X^\ell_{\lambda})-f(X^\ell_{\lambda_0}))|\mu_\mS(\dif\ell)=0.\label{EC99}
\end{align}
The desired continuity follows by noticing
$$
\int_{\mR^d}|\rho_{\lambda_n}(x)-\rho_\lambda(x)|\dif x=\sup_{\|f\|_\infty\leq 1}|\mE(f(X_{\lambda_n})-f(X_{\lambda}))|
$$
and (\ref{EC99}).
\end{proof}

The following theorem about the criterion that a random vector admits a smooth density in the Malliavin calculus can be found in \cite[p.100-103]{Nu}.
\bt\label{Th2}
Assume that $F=(F^1,\cdots, F^d)\in\mD^\infty$ is a smooth functional on $\Omega$ and satisfies that for all $p\geq 2$,
$$
\mE[(\det\Sigma_F)^{-p}]<\infty.
$$
Let $G\in\mD^\infty$ and $\varphi\in C^\infty_b(\mR^d)$. Then for any multi-index $\alpha=(\alpha_1,\cdots,\alpha_m)\in \{1,2,\cdots,d\}^m$,
$$
\mE[\p_\alpha\varphi(F)G]=\mE[\varphi(F)H_{\alpha}(F,G)],
$$
where $\p_\alpha=\p_{\alpha_1}\cdots\p_{\alpha_m}$, and $H_{\alpha}(F,G)$ are recursively defined by
\begin{align}
H_{(i)}(F,G)&:=\sum_{j}D^*\(G(\Sigma^{-1}_F)_{ij}DF^j\),\label{Me1}\\
H_\alpha(F,G)&:=H_{(\alpha_m)}(F,H_{(\alpha_1,\cdots,\alpha_{m-1})}(F,G)).\no
\end{align}
As a consequence, for any $p\geq 1$, there exist $p_1,p_2,p_3>1$ and $n_1,n_2\in\mN$ such that
\begin{align}
\|H_\alpha(F,G)\|_p\leq C\|(\det\Sigma_F)^{-1}\|^{n_1}_{p_1}\|DF\|_{m,p_2}^{n_2}\|G\|_{m,p_3}.\label{Esu}
\end{align}
In particular, the law of $F$ possesses an infinitely differentiable density $\rho\in\cS(\mR^d)$.
\et

\subsection{Malliavin covariance matrix}
Consider the following SDE driven by subordinated Brownian motion $W_{S_t}$:
$$
\dif X_t(x)=b(X_t(x))\dif t+\sigma(X_{t-}(x))\dif W_{S_t}, \ \ X_0(x)=x,
$$
where $b:\mR^d\rightarrow\mR^d$ and $\sigma:\mR^d\rightarrow\mR^d\times\mR^d$ are $C^{\infty}$-functions with bounded partial derivatives of first order.
It is well known that $\{X_t(x), t\geq 0, x\in\mR^d\}$ forms a $C^1$-stochastic flow (cf. \cite{Pr}).
For $0\leq r<t$, let $J_{r,t}(x)$ solve the following linear matrix-valued SDE:
\begin{align}
J_{r,t}(x)=\mI+\int^t_r\nabla b(X_s(x))J_{r,s}(x)\dif s+\sum_{k=1}^d\int_r^t\nabla\sigma_{\cdot k}(X_{s-}(x)) J_{r,s-}(x)\dif W^k_{S_s}.\label{EC9}
\end{align}
Here and below, $\int^t_r$ always means $\int_{(r,t]}$.
We have the following explicit expression about the Malliavin covariance matrix of $X_t(x)$.
\bl
Let $(\Sigma_t(x))_{ij}$ be the Malliavin covariance matrix of $X_t(x)$ defined by (\ref{Def2}). We have
\begin{align}
\Sigma_t(x)=\int^t_0J_{s,t}(x)\sigma(X_{s-}(x))\sigma^*(X_{s-}(x))J^*_{s,t}(x)\dif S_s.\label{Ep4}
\end{align}
\el
\begin{proof}
Fix a path $\ell\in\mS$ and $x\in\mR^d$. Let $X^\ell_t$ solve the following SDE:
\begin{align}
\dif X^\ell_t=b(X^\ell_t)\dif t+\sigma(X^\ell_t)\dif W_{\ell_t}, \ \ X^\ell_0=x,\label{Eq3}
\end{align}
and for $0\leq r<t$, let $J^\ell_{r,t}$ solve the following linear matrix-valued SDE:
\begin{align}
J^\ell_{r,t}=\mI+\int^t_r\nabla b(X^\ell_s)J^\ell_{r,s}\dif s+\sum_{k=1}^d\int_r^t\nabla\sigma_{\cdot k}(X^\ell_{s-}) J^\ell_{r,s-}\dif W^k_{\ell_s}.\label{Eq1}
\end{align}
By definitions (\ref{Def1}) and (\ref{Def2}), it suffices to prove that
\begin{align}
\Sigma^\ell_t=\int^t_0J^\ell_{s,t}\sigma(X^\ell_{s-})\sigma^*(X^\ell_{s-})(J^\ell_{s,t})^*\dif \ell_s,\label{EQ1}
\end{align}
where $(\Sigma^\ell_t)_{ij}:=\<D(X^\ell_t)^i, D(X^\ell_t)^j\>_\mH$.
For any $h\in\mH$, it is by now standard to prove that
\begin{align}
D_h X^\ell_t=\int^t_0\nabla b(X^\ell_s)D_hX^\ell_s\dif s+\sum_k\int^t_0\nabla\sigma_{\cdot k}(X^\ell_{s-})D_hX^\ell_{s-}\dif W^k_{\ell_s}+
\int^t_0\sigma(X^\ell_{s-})\dif h_{\ell_s}.\label{EQ11}
\end{align}
By the formula of constant variations or Fubini's theorem, we have
$$
D_h X^\ell_t=\int^t_0J^\ell_{s,t}\sigma(X^\ell_{s-})\dif h_{\ell_s}.
$$
Let $\{h^n,n\in\mN\}$ be an orthonormal basis of $\mH$. We have
$$
\Sigma^\ell_t=\sum_nD_{h^n}X^\ell_t\otimes D_{h^n}X^\ell_t
=\sum_{n}\int^t_0J^\ell_{s,t}\sigma(X^\ell_{s-})\dif h^n_{\ell_s}\otimes\int^t_0J^\ell_{s,t}\sigma(X^\ell_{s-})\dif h^n_{\ell_s},
$$
where for a vector $x\in\mR^d$, $x\otimes x:=(x_ix_j)_{ij}$.
From this, we then obtain (\ref{EQ1}) by noticing that (cf. \cite[Lemma 4.4]{Zh3})
$$
\sum_n\left(\int^t_0\<f_s,\dif h^n_{\ell_s}\>\right)\left(\int^t_0\<g_s,\dif h^n_{\ell_s}\>\right)=\int^t_0\<f_s, g_s\>\dif \ell_s.
$$
The proof is complete.
\end{proof}

\subsection{Three useful lemmas}
In the remainder of this paper, we make the following assumption on $\nu_S$: For some $\alpha\in(0,2)$,
\begin{align}
0<\varliminf_{\eps\downarrow 0}\eps^{\frac{\alpha}{2}-1}\int^\eps_0u\nu_S(\dif u)
\leq\varlimsup_{\eps\downarrow 0}\eps^{\frac{\alpha}{2}-1}\int^\eps_0u\nu_S(\dif u)<\infty.\tag{\mbox{\bf AS}$^\alpha$}
\end{align}
We use the following filtration:
$$
\sF_t:=\sigma\{W_{S_s},S_s: s\leq t\}.
$$
Clearly, for $t>s$, $W_{S_t}-W_{S_s}$ and $S_t-S_s$ are independent of $\sF_s$. 

We need the following lemmas.
\bl\label{Le10}
Let $f_t:\mR_+\to\mR_+$ be a bounded continuous $\sF_t$-adapted process. Under {\bf(AS)$^\alpha$}, there exist constants $c_0=c_0(\nu_S,\alpha,\|f\|_\infty)\in(0,1)$ and 
$\lambda_0=\lambda_0(\nu_S,\alpha,\|f\|_\infty)>0$ such that for any $\lambda\geq\lambda_0$,
$$
\mE\exp\left\{-\lambda\int^t_0f_s\dif S_s\right\}\leq\left(\mE\exp\left\{-c_0\lambda^{\frac{\alpha}{2}}\int^t_0f_s\dif s\right\}\right)^{\frac{1}{2}}.
$$
\el
\begin{proof}
For $\lambda>0$, set
$$
g^\lambda_s:=\int^\infty_0(1-\mathrm{e}^{-\lambda f_s u})\nu_S(\dif u)
$$
and
$$
M^\lambda_t:=-\lambda\int^t_0f_s\dif S_s+\int^t_0g^\lambda_s\dif s.
$$
Let $N_S(\dif t,\dif u)$ be the Poisson random measure associated with $S_t$, i.e.,
\begin{align}
N_S((0,t]\times U):=\sum_{s\leq t}1_U(S_s-S_{s-}),\ \ U\in\sB(\mR_+).\label{Mu}
\end{align}
Let $\tilde N_S(\dif t,\dif u)$ be the compensated Poisson random measure of $N_S(\dif t,\dif u)$, i.e.,
$$
\tilde N_S(\dif t,\dif u)=N_S(\dif t,\dif u)-\dif t\nu_S(\dif u).
$$
By L\'evy-It\^o's decomposition (cf. \cite{Sa}), we can write
\begin{align}
S_t=t\int^1_0u\nu_S(\dif u)+\int^t_0\!\!\!\int^1_0u \tilde N_S(\dif s,\dif u)+\int^t_0\!\!\!\int^\infty_1u N_S(\dif s,\dif u),\label{EUY}
\end{align}
so that
\begin{align*}
\int^t_0f_s\dif S_s=\int^t_0 f_s\left(\int^1_0u\nu_S(\dif u)\right)\dif s
+\int^t_0\!\!\!\int^1_0f_s u\tilde N_S(\dif s,\dif u)+\int^t_0\!\!\!\int^\infty_1f_s u N_S(\dif s,\dif u).
\end{align*}
By It\^o's formula (cf. \cite{Ap}), we have
\begin{align}
\mathrm{e}^{M^\lambda_t}=1+\int^t_0\!\!\!\int^\infty_0\mathrm{e}^{M^\lambda_{s-}}
[\mathrm{e}^{-\lambda f_su}-1]\tilde N_S(\dif s,\dif u).\label{EQ2}
\end{align}
Since for any $x\geq 0$,
$$
1-\mathrm{e}^{-x}\leq 1\wedge x,
$$
we have
$$
M^\lambda_t\leq \int^t_0g^\lambda_s\dif s
\leq t\int^\infty_0(1\wedge(\lambda \|f\|_\infty |u|))\nu_S(\dif u).
$$
Hence, by (\ref{EQ2}) we have
$$
\mE \mathrm{e}^{M^\lambda_t}=1,
$$
and by H\"order's inequality,
\begin{align}
\mE\exp\left\{-\frac{\lambda}{2}\int^t_0f_s\dif S_s\right\}&=\mE \left( \e^{M_t^{\lambda}/2}\exp\left\{-\frac{1}{2}\int^t_0g^{\lambda}_s\dif s\right\}\right)
\leq\left(\mE \exp\left\{-\int^t_0g^{\lambda}_s\dif s\right\}\right)^{\frac{1}{2}}.\label{EK6}
\end{align}
On the other hand, since for any $\kappa\in(0,1)$ and $x\leq-\log k$,
$$
1-\mathrm{e}^{-x}\geq \kappa x,
$$
letting $\kappa=\frac{1}{\e}$, we have
\begin{align*}
g^{\lambda}_s&\geq\int^{\frac{1}{\lambda\|f\|_\infty}}_0(1-\mathrm{e}^{-\lambda f_s u})\nu_S(\dif u)
\geq\frac{1}{\e}\int^{\frac{1}{\lambda \|f\|_\infty}}_0(\lambda f_s u)\nu_S(\dif u)=\frac{\lambda f_s}{\e}\int^{\frac{1}{\lambda \|f\|_\infty}}_0u\nu_S(\dif u).
\end{align*}
Substituting this into (\ref{EK6}) and using {\bf(AS$^\alpha$)}, we obtain the desired estimate.
\end{proof}
\bl
Under {\bf(AS$^\alpha$)}, there exist constants $C_0, C_1\geq 1$ such that for all $\delta\in(0,1)$,
\begin{align}
\int_{|z|\geq\delta}\nu_L(\dif z)\leq C_0\delta^{-\alpha},\ \ 
\int_{|z|\leq\delta}|z|^2\nu_L(\dif z)&\leq C_1\delta^{2-\alpha}.\label{EUT0}
\end{align}
\el
\begin{proof}
First of all, by {\bf(AS$^\alpha$)}, there exists a constant $C\geq 1$ such that for all $\eps\in(0,1)$,
\begin{align}
\int^\eps_0u\nu_S(\dif u)\leq C\eps^{1-\frac{\alpha}{2}}.\label{EC8}
\end{align}
For any $\eps\in(0,1)$, letting $m:=[\log_2(1/\eps)]+1$, we have
\begin{align*}
\nu_S([\eps,1])&\leq\sum_{k=0}^m\int^{2^{k+1}\eps}_{2^k\eps}\nu_S(\dif u)
\leq\sum_{k=0}^m(2^k\eps)^{-1}\int^{2^{k+1}\eps}_{2^k\eps}u\nu_S(\dif u)\\
&\leq C\sum_{k=0}^m(2^k\eps)^{-1}(2^{k+1}\eps)^{1-\frac{\alpha}{2}}=C2^{1-\frac{\alpha}{2}}\sum_{k=0}^m (2^k\eps)^{-\frac{\alpha}{2}},
\end{align*}
which implies that for some $C\geq 1$,
\begin{align}
\nu_S([\eps,\infty))\leq C\eps^{-\frac{\alpha}{2}}.\label{EC99}
\end{align}
By (\ref{Levy}) and the change of variables, we have
\begin{align*}
\int_{|z|\geq\delta}\nu_L(\dif z)&=\int^\infty_0\left(\int_{|z|\geq\delta}(2\pi u)^{-\frac{d}{2}}\e^{-\frac{|z|^2}{2u}}\dif z\right)\nu_S(\dif u)\\
&\leq\nu_S([\delta^2,\infty))+\pi^{-\frac{d}{2}}\int^{\delta^2}_0\left(\int_{\sqrt{2u}|x|\geq\delta}\e^{-|x|^2}\dif x\right)\nu_S(\dif u)\\
&\stackrel{(\ref{EC99})}{\leq} C\delta^{-\alpha}+\pi^{-\frac{d}{2}}
\int^{\delta^2}_0u\left(\int_{\sqrt{2u}|x|\geq\delta}\frac{2|x|^2}{\delta^2}\e^{-|x|^2}\dif x\right)\nu_S(\dif u)\\
&\leq C\delta^{-\alpha}+2\pi^{-\frac{d}{2}}\delta^{-2}\int^{\delta^2}_0u\nu_S(\dif u)\left(\int_{\mR^d}|x|^2\e^{-|x|^2}\dif x\right)
\stackrel{(\ref{EC8})}{\leq} C_0\delta^{-\alpha}
\end{align*}
and
\begin{align*}
\int_{|z|\leq\delta}|z|^2\nu_L(\dif z)&=\int^\infty_0\left(\int_{|z|\leq\delta}(2\pi u)^{-\frac{d}{2}}|z|^2\e^{-\frac{|z|^2}{2u}}\dif z\right)\nu_S(\dif u)\\
&=2\pi^{-\frac{d}{2}}\int^\infty_0 u\left(\int_{\sqrt{2u}|x|\leq\delta}|x|^2\e^{-|x|^2}\dif x\right)\nu_S(\dif u)\\
&\leq C_d\left(\int^{\delta^2}_0 u\nu_S(\dif u)+\delta^2\int^\infty_{\delta^2}\nu_S(\dif u)\right)\stackrel{(\ref{EC8})(\ref{EC99})}{\leq} C_1\delta^{2-\alpha}.
\end{align*}
The proof is complete.
\end{proof}
\bl\label{Le26}
Under {\bf(AS$^\alpha$)}, we have
$$
\mP\Big(\omega: \{s: \Delta S_s(\omega)>0\} \mbox{ is dense in }(0,\infty)\Big)=1.
$$
\el
\begin{proof}
Define a stopping time $\tau:=\inf\{t>0: S_t>0\}=\inf\{t>0: \Delta S_t>0\}$. As in the proof of \cite[Lemma 2.1]{Zh2}, it suffices to prove that
$$
\mP(\tau=0)=1.
$$
Let $N_S(\dif s,\dif u)$ be defined by (\ref{Mu}). For any $\eps\in(0,1)$, we have
\begin{align*}
\eps\geq \mE\left(\Delta S_\tau1_{\Delta S_\tau\leq\eps}\right)
&=\mE\left(\sum_{0<t\leq\tau}\Delta S_t 1_{\Delta S_t\leq\eps}\right)=\mE\left(\int^{\tau}_0\!\!\!\int^\eps_0 u N_S(\dif s,\dif u)\right)\\
&=\mE\left(\int^{\tau}_0\!\!\!\int^\eps_0u\nu_S(\dif z)\dif s\right)
=\left(\int^\eps_0 u\nu_S(\dif u)\right)\mE\tau,
\end{align*}
which, together with {\bf(AS$^\alpha$)} and letting $\eps\to 0$, implies that
$$
\mE\tau=0\Rightarrow \mP(\tau=0)=1.
$$
The proof is complete.
\end{proof}

\section{Proof of Theorem \ref{Th00}}

In order to prove Theorem \ref{Th00}, it suffices to check the conditions of Theorem \ref{Th11}, and (\ref{RT1}) follows by It\^o's formula. First of all, we have
\bl
Assume that $b,\sigma\in C^\infty$ has bounded partial derivatives of first order. For each $\ell\in\mS$, $p\geq 2$ and all $(t,x)\in(0,\infty)\times\mR^d$, we have
\begin{align}
X^\ell_t(x)\in\mD^{1,p}(\mW;\mR^d).\label{CV1}
\end{align}
Moreover, for fixed $t>0$ and $x\in\mR^d$, there exists a $\mu_\mS$-null set $\mS_0$ such that for all $\ell\notin \mS_0$,
\begin{align}
\lim_{(s,y)\to(t,x)}\|X^\ell_s(y)-X^\ell_t(x)\|_{1,p}=0.\label{CV2}
\end{align}
In particular, {\bf (H1)} of Theorem \ref{Th11} holds.
\el
\begin{proof}
First of all, by equations (\ref{EQ1}) and (\ref{EQ11}), since $b,\sigma\in C^1_b$, it is easy to prove that for any $p\geq 2$ and $T>0$, $x\in\mR^d$,
\begin{align*}
\mE\left(\sup_{t\in[0,T]}(|X^\ell_t(x)|^p+\|DX^\ell_t(x)\|^p_\mH)\right)\leq C(T,p,\ell)(1+|x|^p),
\end{align*} 
which gives (\ref{CV1}).
Let us prove (\ref{CV2}). 
For fixed $x,y\in\mR^d$, set
$$
Z_t:=X^\ell_t(x)-X^\ell_t(y).
$$
By equation (\ref{EQ11}), we have
\begin{align*}
\<DZ_t,h\>_\mH&=x-y+\int^t_0\nabla b(X_s(x)) \<DZ_s,h\>_\mH\dif s
+\sum_k\int^t_0\nabla\sigma_{\cdot k}(X_{s-}(x))\<DZ_{s-},h\>_\mH\dif W^k_{\ell_s}\\
&+\int^t_0[\nabla b(X_s(x))-\nabla b(X_s(y))] D_hX_s(y)\dif s+\int^t_0[\sigma(X_{s-}(x))-\sigma(X_{s-}(y))]\dif h_{\ell_s}\\
&+\sum_k\int^t_0[\nabla\sigma_{\cdot k}(X_{s-}(x))-\nabla\sigma_{\cdot k}(X_{s-}(y))] D_hX_{s-}(y)\dif W^k_{\ell_s}.
\end{align*}
Let $\{h^n,n\in\mN\}$ be an orthonormal basis of $\mH$. Since
$$
\|DZ_t\|_\mH^2=\sum_{n}|\<DZ_t,h_n\>_\mH|^2,
$$
by Burkholder's inequality (cf. \cite{Ku} or \cite[Lemma 2.3]{Zh1}), we have
\begin{align*}
\mE\|DZ_t\|^p_\mH&\leq C|x-y|^p+Ct^{p-1}\int^t_0\mE\|DZ_s\|^p_\mH\dif s+C\mE\left(\int^t_0\|DZ_{s-}\|^2_\mH\dif \ell_s\right)^{p/2}\\
&+Ct^{p-1}\int^t_0\mE[|\nabla b(X_s(x))-\nabla b(X_s(y))|^p \|DX_s(y)\|_\mH^p]\dif s+C\mE\left(\int^t_0|Z_{s-}|^2\dif \ell_s\right)^{p/2}\\
&+C\mE\left(\int^t_0|\nabla\sigma(X_{s-}(x))-\nabla\sigma(X_{s-}(y))|^2 \|DX_{s-}(y)\|_\mH^2\dif \ell_s\right)^{p/2}\\
&\leq C|x-y|^p+Ct^{p-1}\int^t_0\mE\|DZ_s\|^p_\mH\dif s+C\ell^{p/2-1}_t\int^t_0\mE\|DZ_{s-}\|^p_\mH\dif \ell_s+I_t(x,y),
\end{align*}
where $C$ only depends on $p$, and $I_t(x,y)$ denotes the remaining terms. By Gronwall's inequality, we derive that
$$
\mE\|DZ_t\|^p_\mH\leq C^\ell_t(|x-y|^p+I_t(x,y)).
$$
On the other hand, we similarly have
$$
\mE|X^\ell_t(x)-X^\ell_t(y)|^p\leq C^\ell_t |x-y|^p
$$
and
$$
\mE\|X^\ell_t(x)-X^\ell_s(x)\|^p_{1,p}\leq C(|t-s|^p+|\ell_t-\ell_s|^{p/2})(1+|x|^p).
$$
Combining the above estimates and noticing that for fixed $t$ and $\mu_\mS$-almost all $\ell$, 
$s\mapsto\ell_s$ is continuous at $t$, we thus obtain the desired continuity (\ref{CV2}).
\end{proof}

The remaining task is to verify {\bf (H2)} of Theorem \ref{Th11} under (\ref{ForH}), i.e., 
to prove the invertibility of the  Malliavin covariance matrix $\Sigma_t(x)$ given in (\ref{Ep4}). 

Let $N(\dif t,\dif z)$ be the Poisson random measure associate to $L_t=W_{S_t}$, i.e.,
$$
N((0,t]\times U)=\sum_{s\leq t}1_U(\Delta L_s),\ \ U\in\sB(\mR^d_0).
$$
The compensated Poisson random measure of $N(\dif t,\dif z)$ is defined by
$$
\tilde N(\dif t,\dif z):=N(\dif t,\dif z)-\dif t\nu_L(\dif z).
$$
Choose $\delta>0$ being small and fixed so that
\begin{align}
2\delta\|\nabla\sigma\|_\infty\leq 1.\label{Del}
\end{align}
Define
\begin{align}
L^\delta_t:=\int^t_0\!\!\!\int_{|z|\leq\delta}z\tilde N(\dif s,\dif z),\ \ \hat L^\delta_t:=\int^t_0\!\!\!\int_{|z|>\delta}z N(\dif s,\dif z).\label{NB4}
\end{align} 
Then $L^\delta_t$ and $\hat L^\delta_t$ are two independent L\'evy processes with L\'evy measures $1_{\{|z|\leq\delta\}}\nu_L(\dif z)$ and
$1_{\{|z|>\delta\}}\nu_L(\dif z)$ respectively, and by L\'evy-It\^o's decomposition, we have
\begin{align}
L_t=L^\delta_t+\hat L^\delta_t.\label{GH7}
\end{align}
Fix $T>0$ and a c\`adl\`ag function $\hbar:(0,\infty)\to\mR^d$ with 
\begin{align}
\Delta\hbar_T=0,\ \ |\Delta \hbar_s|=0\mbox{ or }>\delta,\ \ \forall s\geq 0.\label{GH3}
\end{align}
Notice that for almost all $\omega$, the path $t\mapsto \hat L^\delta_t(\omega)$ has this property.
Let $X^\hbar_t$ solve the following SDE:
$$
X^\hbar_t=x+\int^t_0b(X^\hbar_s)\dif s+\int^t_0\!\!\!\int_{|z|\leq\delta}\sigma(X^\hbar_{s-})z\tilde N(\dif s,\dif z)+\int^t_0\sigma(X^\hbar_{s-})\dif \hbar_s,
$$
and $J^\hbar_{r,t}$ solve the following linear matrix-valued SDE:
\begin{align*}
J^\hbar_{r,t}&=\mI+\int^t_r\nabla b(X^\hbar_s)J^\hbar_{r,s}\dif s+\int_r^t\!\!\!\int_{|z|\leq\delta}\nabla\sigma(X^\hbar_{s-})z J^\hbar_{r,s-}\tilde N(\dif s,\dif z)
+\sum_{k=1}^d\int_r^t\nabla\sigma_{\cdot k}(X^\hbar_{s-}) J^\hbar_{r,s-}\dif \hbar^k_s.
\end{align*}
By (\ref{GH7}), it is easy to see that
\begin{align}
X_t=X^\hbar_t |_{\hbar=\hat L^\delta_\cdot},\ \ J_{r,t}=J^\hbar_{r,t}|_{\hbar=\hat L^\delta_\cdot}.\label{GH2}
\end{align}
Define
\begin{align}
\lambda:=\sup\{t<T: |\Delta\hbar_t|>\delta\}.\label{Tau}
\end{align}
By (\ref{GH3}), we have for all $t\in[\lambda,T]$,
$$
J^\hbar_{\lambda,t}=\mI+\int^t_\lambda\nabla b(X^\hbar_s)J^\hbar_{\lambda,s}\dif s+\int_\lambda^t\!\!\!\int_{|z|\leq\delta}\nabla\sigma(X^\hbar_{s-})z J^\hbar_{\lambda,s-}\tilde N(\dif s,\dif z).
$$
Define
$$
Q(x,z):=(\mI+\nabla\sigma(x) z)^{-1}-\mI,\ \ |z|\leq\delta.
$$
By (\ref{Del}), this is well defined and
\begin{align}
|Q(x,z)|\leq C_d\|\nabla\sigma\|_\infty|z|,\ \ |z|\leq\delta.\label{Es2}
\end{align}
We have
\bl\label{Le1}
For each $t\in[\lambda,T]$, the matrix $J^\hbar_{\lambda,t}$ is invertible and $K^\hbar_{\lambda,t}:=(J^\hbar_{\lambda,t})^{-1}$ solves
\begin{align*}
K^\hbar_{\lambda,t}&=\mI-\int^t_\lambda K^\hbar_{\lambda,s}\nabla b(X^\hbar_s)\dif s
+\int_\lambda^t\!\!\!\int_{|z|\leq\delta}K^\hbar_{\lambda,s-}Q(X^\hbar_{s-},z)\tilde N(\dif s,\dif z)\no\\
&\qquad-\int_\lambda^t\!\!\!\int_{|z|\leq\delta}K^\hbar_{\lambda,s-}Q(X^\hbar_{s-},z)\nabla\sigma(X^\hbar_{s-})z\nu_L(\dif z)\dif s.
\end{align*}
Moreover, for any $p\geq 2$, we have
\begin{align}
\sup_{t\in[\lambda,T]}\mE|K^\hbar_{\lambda,t}|^p<\infty.\label{Es1}
\end{align}
\el
\begin{proof}
For simplicity of notations, we drop the superscript ``$\hbar$''.
By It\^o's formula, we have
\begin{align*}
K_{\lambda,t}J_{\lambda,t}&=\mI+\int^t_\lambda K_{\lambda,s-}\dif J_{\lambda,s}+\int^t_\lambda \dif K_{\lambda,s} J_{\lambda,s-}+\sum_{\lambda<s\leq t}\Delta K_{\lambda,s}\Delta J_{\lambda,s}\\
&=\mI+\int^t_\lambda K_{\lambda,s-}\nabla b(X_{s-}) J_{\lambda,s-}\dif s+\int_\lambda^t\!\!\!\int_{|z|\leq\delta}K_{\lambda,s-}\nabla \sigma(X_{s-})z J_{\lambda,s-}\tilde N(\dif s,\dif z)\\
&\quad-\int^t_\lambda K_{\lambda,s-}\nabla b(X_{s-}) J_{\lambda,s-}\dif s+\int_\lambda^t\!\!\!\int_{|z|\leq\delta}K_{\lambda,s-}Q(X_{s-},z)J_{\lambda,s-}\tilde N(\dif s,\dif z)\\
&\quad-\int_\lambda^t\!\!\!\int_{|z|\leq\delta}K_{\lambda,s-}Q(X_{s-},z)\nabla\sigma(X_{s-})zJ_{\lambda,s-}\nu_L(\dif z)\dif r\\
&\quad+\int_\lambda^t\!\!\!\int_{|z|\leq\delta}K_{\lambda,s-}Q(X_{s-},z)\nabla\sigma(X_{r-})zJ_{\lambda,s-} N(\dif s,\dif z).
\end{align*}
By the definition of $Q(x,z)$, one sees that
$$
Q(x,z)+\nabla\sigma(x)z=-Q(x,z)\nabla\sigma(x)z.
$$
Hence,
$$
K_{\lambda,t}J_{\lambda,t}=\mI.
$$
Estimate (\ref{Es1}) is standard by (\ref{Es2}) and Gronwall's inequality. The proof is complete.
\end{proof}

Let $V:\mR^d\to\mR^d\times\mR^d$ be a matrix-valued smooth function. Let us first introduce two functions:
\begin{align}
H_V(x,z)&:=V(x+\sigma(x)z)-V(x)+Q(x,z)V(x+\sigma(x)z),\label{Re12}\\
G_V(x,z)&:=H_V(x,z)+\nabla\sigma(x)z\cdot V(x)-\sigma(x)z\cdot\nabla V(x).\label{Re22}
\end{align}
Define
\begin{align}
M_V^\hbar(t)&:=\int_\lambda^t\!\!\!\int_{|z|\leq\delta}K^\hbar_{\lambda,s-}H_V(X^\hbar_{s-},z)\tilde N(\dif s,\dif z),\label{Re32}\\
R_V^\hbar(t)&:=\int_\lambda^t\!\!\!\int_{|z|\leq\delta}K^\hbar_{\lambda,s}G_V(X^\hbar_{s-},z)\nu_L(\dif z)\dif s.\label{Re42}
\end{align}
\bl\label{Le32}
Let $[b,V]:=b\cdot\nabla V-V\cdot\nabla b$. We have $\mP(\Omega^V_1)=1$, where
\begin{align*}
\Omega^V_1:=\left\{\omega: K^\hbar_{\lambda,t}V(X^\hbar_t)=V(X^\hbar_\lambda)+\int^t_\lambda K^\hbar_{\lambda,s}[b,V](X^\hbar_s)\dif s+M^\hbar_V(t)+R^\hbar_V(t),
\forall t\in[\lambda,T]\right\}.
\end{align*}
\el
\begin{proof}
For simplicity of notations, we drop the superscript ``$\hbar$''. By It\^o's formula, we have
\begin{align*}
V(X_t)&=V(X_\lambda)+\int^t_\lambda (b\cdot\nabla) V(X_s)\dif s
+\int_\lambda^t\!\!\!\int_{|z|\leq\delta}[V(X_{s-}+\sigma(X_{s-})z)-V(X_{s-})]\tilde N(\dif s,\dif z)\\
&\quad+\int_\lambda^t\!\!\!\int_{|z|\leq\delta}[V(X_{s-}+\sigma(X_{s-})z)-V(X_{s-})-\sigma(X_{s-})z\cdot\nabla V(X_{s-})]\nu_L(\dif z)\dif s.
\end{align*}
By It\^o's product formula and Lemma \ref{Le1}, we have
\begin{align*}
K_{\lambda,t}V(X_t)&=V(X_\lambda)+\int_\lambda^t\dif K_{\lambda, s}  V(X_{s-})+\int_\lambda^tK_{\lambda, s-}\dif V(X_s)
+\sum_{\lambda<s\leq t}\Delta K_{\lambda,s}\Delta V(X_s)\\
&=V(X_\lambda)+\int_\lambda^tK_{\lambda, s-}[b,V](X_s)\dif s
+\int_\lambda^t\!\!\!\int_{|z|\leq\delta}K_{\lambda,s-}Q(X_{\lambda,s-},z)V(X_{s-})\tilde N(\dif s,\dif z)\\
&\quad-\int_\lambda^t\!\!\!\int_{|z|\leq\delta}K_{\lambda,s-}Q(X_{\lambda,s-},z)\nabla\sigma(X_{\lambda,s-})zV(X_{s-})\nu_L(\dif z)\dif s\\
&\quad+\int_\lambda^t\!\!\!\int_{|z|\leq\delta}K_{\lambda, s-}[V(X_{s-}+\sigma(X_{s-})z)-V(X_{s-})]\tilde N(\dif s,\dif z)\\
&\quad+\int_\lambda^t\!\!\!\int_{|z|\leq\delta}K_{\lambda, s-}[V(X_{s-}+\sigma(X_{s-})z)-V(X_{s-})-\sigma(X_{s-})z\cdot\nabla V(X_{s-})]\nu_L(\dif z)\dif s\\
&\quad+\int_\lambda^t\!\!\!\int_{|z|\leq\delta}K_{\lambda, s-}Q(X_{s-},z)[V(X_{s-}+\sigma(X_{s-})z)-V(X_{s-})]N(\dif s,\dif z).
\end{align*}
Noticing that
\begin{align*}
&V(x+\sigma(x)z)-V(x)-\sigma(x)z\cdot\nabla V(x)+Q(x,z)[V(x+\sigma(x)z)-V(x)-\nabla\sigma(x)z V(x)]\\
&\quad=V(x+\sigma(x)z)-V(x)+Q(x,z)V(x+\sigma(x)z)+\nabla\sigma(x)z\cdot V(x)-\sigma(x)z\cdot\nabla V(x)\\
&\quad=H_V(x,z)+\nabla\sigma(x)z\cdot V(x)-\sigma(x)z\cdot\nabla V(x)=G_V(x,z),
\end{align*}
where $H_V$ and $G_V$ are defined by (\ref{Re12}) and (\ref{Re22}), we obtain
\begin{align}
K_{\lambda,t}V(X_t)&=V(X_\lambda)+\int_\lambda^tK_{\lambda, s-}[b,V](X_s)\dif s
+\int_\lambda^t\!\!\!\int_{|z|\leq\delta}K_{\lambda,s-}H_V(X_{s-},z)\tilde N(\dif s,\dif z)\no\\
&\qquad+\int_\lambda^t\!\!\!\int_{|z|\leq\delta}K_{\lambda,s-}G_V(X_{s-},z)\nu_L(\dif z)\dif s,\label{IT0}
\end{align}
which gives the desired result.
\end{proof}

We have
\bl\label{Le33}
There exists a sequence of numbers $\eps_m\in(0,\delta)\to 0$ such that $\mP(\Omega^V_2)=1$, where
\begin{align*}
\Omega^V_2:=\Bigg\{\omega: 
M^\hbar_V(t)=\lim_{m\to\infty}\int_\lambda^t\!\!\!\int_{\eps_m<|z|\leq\delta}
K^\hbar_{\lambda,s-}H_V(X^\hbar_{s-},z)\tilde N(\dif s,\dif z)\mbox{ uniformly in $t\in[\lambda,T]$}\Bigg\}.
\end{align*}
\el
\begin{proof}
Since $b$ and $\sigma$ are linear growth, by using a standard stopping time technique, we may assume that $\sigma, V$ and $\nabla V$ are bounded so that
\begin{align}
|H_V(x,z)|\leq \|\nabla V\|_\infty\|\sigma\|_\infty |z|+C_d\|\nabla\sigma\|_\infty\|V\|_\infty |z|. \label{GH9}
\end{align}
By Doob's maximal inequality, we have
\begin{align*}
&\mE\left(\sup_{t\in[\lambda,T]}\left|M^\hbar_V(t)-\int_\lambda^t\!\!\!\int_{\eps<|z|\leq\delta}K^\hbar_{\lambda,s-}H_V(X^\hbar_{s-},z)\tilde N(\dif s,\dif z)\right|^2\right)\\
&\qquad=\mE\left(\sup_{t\in[\lambda,T]}\left|\int_\lambda^t\!\!\!\int_{|z|\leq\eps}K^\hbar_{\lambda,s-}H_V(X^\hbar_{s-},z)\tilde N(\dif s,\dif z)\right|^2\right)\\
&\qquad\leq 4\mE\left|\int_\lambda^T\!\!\!\int_{|z|\leq \eps}K^\hbar_{\lambda,s-}H_V(X^\hbar_{s-},z)\tilde N(\dif s,\dif z)\right|^2\\
&\qquad\leq4\mE\left(\int^T_\lambda\!\int_{|z|\leq\eps}|K^\hbar_{\lambda,s}H_V(X^\hbar_{s-},z)|^2\nu_L(\dif z)\dif s\right),
\end{align*}
which tends to zero as $\eps\downarrow 0$ by (\ref{Es1}) and (\ref{GH9}).
We complete the proof by choosing a suitable sequence $\eps_n\downarrow 0.$
\end{proof}

By \cite[p.68, Theorem 23]{Pr}, we also have
\bl\label{Le34}
For $n\in\mN$, let $t_k:=(kT)/n\vee\lambda$. There exists a subsequence $n_m\to\infty$ such that  $\mP(\Omega^V_3)=1$, where
\begin{align*}
\Omega^V_3&:=\Bigg\{\omega:\  
\lim_{n_m\to\infty}\sum_{k=0}^{n_m-1}(M^\hbar_V(t_{k+1}\wedge t)-M^\hbar_V(t_k\wedge t))_{ij}(M^\hbar_V(t_{k+1}\wedge t)-M^\hbar_V(t_k\wedge t))_{i'j'}\\
&\qquad\quad=\sum_{k,k'}\int^t_\lambda\!\!\!\int_{|z|\leq\delta}(K^\hbar_{\lambda,s-})_{ik}(H_V(X^\hbar_{s-},z))_{kj}
(K^\hbar_{\lambda,s-})_{i'k'}(H_V(X^\hbar_{s-},z))_{k'j'}N(\dif s,\dif z)\\
&\qquad\qquad\mbox{ uniformly in $t\in[\lambda,T]$},\ \forall i,j,i',j'=1,\cdots,d\Bigg\}.
\end{align*}
\el

Now we are in a position to prove the following main result of this section.
\bl
Under (\ref{ForH}), for each $T>0$ and $x\in\mR^d$, the Malliavin covariance matrix $\Sigma_T(x)$ defined by (\ref{Ep4}) is invertible almost surely.
\el
\begin{proof} Below we drop the variable ``$x$'', and divide the proof into four steps.
\\
\\
{\bf (1)} We first prove that 
\begin{align}
t\mapsto J_{t,T}\mbox{ has a c\`adl\`ag modification on $[0,T)$.}\label{GH1}
\end{align}
Let $\tau_0:=0$ and define a sequence of stopping times by
$$
\tau_n:=\{t>\tau_{n-1}: |\Delta L_t|>\delta\},\ \ n\in\mN.
$$
Recalling equation (\ref{EC9}) and by (\ref{GH7}) and (\ref{NB4}), we have
\begin{align}
J_{r,t}&=\mI+\int^t_r\nabla b(X_s)J_{r,s}\dif s+\int_r^t\!\int_{|z|\leq\delta}\nabla\sigma(X_{s-})z J_{r,s-}\tilde N(\dif s,\dif z)\no\\
&\qquad+\int_r^t\!\int_{|z|>\delta}\nabla\sigma(X_{s-})z J_{r,s-} N(\dif s,\dif z).\label{Es4}
\end{align}
Thus, by (\ref{Del}) and as in Lemma \ref{Le1}, we have
$$
J_{\tau_n, t}\mbox{ is invertible for any $t\in[\tau_n,\tau_{n+1})$},
$$
and the inverse $K_{\tau_n,t}:=J^{-1}_{\tau_n, t}$ solves
\begin{align*}
K_{\tau_n,t}&=\mI-\int^t_{\tau_n}K_{\tau_n,s}\nabla b(X_s)\dif s
+\int_{\tau_n}^t\!\int_{|z|\leq\delta}K_{\tau_n,s-}Q(X_{s-},z)\tilde N(\dif s,\dif z)\no\\
&\qquad-\int_{\tau_n}^t\!\int_{|z|\leq\delta}K_{\tau_n,s-}Q(X_{s-},z)\nabla\sigma(X_{s-})z\nu_L(\dif z)\dif s.
\end{align*}
In particular,
\begin{align}
\mbox{$t\mapsto K_{\tau_n,t}$ is c\`adl\`ag on $[\tau_n,\tau_{n+1})$}.\label{GH}
\end{align}
Since for each $t\in[\tau_n,\tau_{n+1})$, by the uniqueness of solutions to equation (\ref{Es4}),
$$
J_{\tau_n,T}=J_{t,T}J_{\tau_n,t}\ \ a.s.,
$$
we have
$$
J_{t,T}=J_{\tau_n,T}J^{-1}_{\tau_n,t}=J_{\tau_n,T}K_{\tau_n,t}\ \ a.s.,
$$
which together with (\ref{GH}) gives (\ref{GH1}).
\\
\\
{\bf (2)} In order to prove the invertibility of $\Sigma_T$, by (\ref{Ep4}) we only need to prove that
$$
\mP\left(\omega: u\Sigma_T u^*=\int^T_0 |uJ_{t,T}\sigma(X_{t-})|^2\dif S_t=0,\ \exists u\in\mS^{d-1}\right)=0.
$$
Notice that
$$
\int^T_0 |uJ_{t,T}\sigma(X_{t-})|^2\dif S_t=\sum_{0<t\leq T}|uJ_{t,T}\sigma(X_{t-})|^2\Delta S_t.
$$
If we set
$$
\sI:=\{t: \Delta S_t>0\},
$$
then it suffices to prove that
$$
\mP\left(\omega: uJ_{t,T}\sigma(X_{t-})=0,\ \ \forall t\in \sI\cap[0,T],\ \exists u\in\mS^{d-1}\right)=0,
$$
which, by Lemma \ref{Le26}, (\ref{GH1}) and the right continuity of $t\mapsto X_t$,  is equivalent to
$$
\mP\left(\omega: uJ_{t,T}\sigma(X_t)=0,\ \ \forall t\in [0,T),\ \exists u\in\mS^{d-1}\right)=0.
$$
Furthermore, by (\ref{GH2}) and taking regular conditional probability with respect to the large jump $\hat L^\delta_\cdot$, 
it is enough to prove that for each c\`adl\`ag $\hbar$ satisfying (\ref{GH3}),
$$
\mP\left(\omega: uJ^\hbar_{t,T}\sigma(X^\hbar_t)=0,\ \ \forall t\in [0,T),\ \exists u\in\mS^{d-1}\right)=0.
$$
Let $\lambda$ be defined by (\ref{Tau}). By Lemma \ref{Le1}, we have for all $t\in[\lambda,T]$,
$$
J^\hbar_{\lambda,T}=J^\hbar_{t, T}J^\hbar_{\lambda,t}\Rightarrow J^\hbar_{t,T}=J^\hbar_{\lambda, T} K^\hbar_{\lambda,t}.
$$
Thus, it reduces to prove
\begin{align}
\mP\left(\omega: uK^\hbar_{\lambda,t}\sigma(X^\hbar_t)=0,\ \ \forall t\in [\lambda,T),\ \exists u\in\mS^{d-1}\right)=0.\label{GH5}
\end{align}
\\
{\bf (3)} Let $B_n$ be defined by (\ref{EUT1}), and 
let $\Omega^{B_n}_1,\Omega^{B_n}_2$ and $\Omega^{B_n}_3$ be defined as in Lemmas \ref{Le32}, \ref{Le33} and \ref{Le34} respectively. Set 
$$
\tilde\Omega:=\cap_{n\in\mN}(\Omega^{B_n}_1\cap\Omega^{B_n}_2\cap\Omega^{B_n}_3).
$$
We want to prove that for each $\omega\in\tilde\Omega$,
\begin{align}
uK^\hbar_{\lambda,t}(\omega)\sigma(X^\hbar_t(\omega))=0,\ \ \forall t\in [\lambda,T)\Rightarrow u=0.\label{RT2}
\end{align}
If this is proven, then by Lemmas \ref{Le32}, \ref{Le33} and \ref{Le34}, we immediately obtain (\ref{GH5}).

Below, for simplicity of notations, we drop ``$\omega$''. Now suppose that
$$
uK^\hbar_{\lambda,t}\sigma(X^\hbar_t)=uK^\hbar_{\lambda,t}B(X^\hbar_t)=0,\ \ \forall t\in [\lambda,T).
$$
We shall use the induction to prove that for any $n\in\mN$,
\begin{align}
uK^\hbar_{\lambda,t}B_n(X^\hbar_t)=0,\ \ \forall t\in [\lambda,T).\label{Eq5}
\end{align}
If this is proven, then by letting $t\downarrow\lambda$, we obtain
$$
uB_n(X^\hbar_\lambda)=0,\ \ \forall n\in\mN,
$$
which implies $u=0$ by (\ref{ForH}), and so (\ref{RT2}).
\\
\\
{\bf (4)} Suppose now that (\ref{Eq5}) holds for some $n\in\mN$. In view of $\omega\in\Omega^{B_n}_1$, by Lemma \ref{Le32} and recalling $B_{n+1}=[b,B_n]$, we have
$$
0=uK^\hbar_{\lambda,t}B_n(X^\hbar_t)=uB_n(X^\hbar_\lambda)+\int^t_\lambda uK^\hbar_{\lambda,s}B_{n+1}(X^\hbar_s)\dif s
+uM^\hbar_{B_n}(t)+uR^\hbar_{B_n}(t),\ \ \forall t\in [\lambda,T).
$$
From this, we have
$$
uB_n(X^\hbar_\lambda)=0
$$
and
\begin{align}
\int^t_\lambda uK^\hbar_{\lambda,s}B_{n+1}(X^\hbar_s)\dif s+uM^\hbar_{B_n}(t)+uR^\hbar_{B_n}(t)=0,\ \ \forall t\in [\lambda,T).\label{Eq4}
\end{align}
Since $t\mapsto \int^t_\lambda uK^\hbar_{\lambda,s}B_{n+1}(X^\hbar_s)\dif s+uR^\hbar_{B_n}(t)$ is absolutely continuous with respect to the Lebesgue measure, 
in view of $\omega\in\Omega^{B_n}_3$, by Lemma \ref{Le34} and (\ref{Eq4}), we have
\begin{align*}
\int_\lambda^t\!\!\!\int_{|z|\leq\delta} |uK^\hbar_{\lambda,s-}H_{B_n}(X^\hbar_{s-},z)|^2N(\dif s,\dif z)=0,\ \ \forall t\in [\lambda,T),
\end{align*}
which implies that for any $\eps\in(0,\delta)$,
\begin{align}
\int_\lambda^t\!\!\!\int_{\eps<|z|\leq\delta} uK^\hbar_{\lambda,s-}H_{B_n}(X^\hbar_{s-},z)N(\dif s,\dif z)=0,\ \ \forall t\in [\lambda,T).\label{Eq33}
\end{align}
On the other hand, since $\nu_L$ is symmetric, by (\ref{Re42}) and (\ref{Re22}), we have for any $t\in[\lambda,T)$,
\begin{align*}
uR^\hbar_{B_n}(t)=\lim_{\eps\downarrow 0}\int_\lambda^t\!\!\!\int_{\eps<|z|\leq\delta}uK^\hbar_{\lambda,s-}H_{B_n}(X^\hbar_{s-},z)\nu_L(\dif z)\dif s,
\end{align*}
and in view of  $\omega\in\Omega^{B_n}_2$, by Lemma \ref{Le33},
$$
uM^\hbar_{B_n}(t)=\lim_{m\to\infty}\int_\lambda^t\!\!\!\int_{\eps_m<|z|\leq\delta}uK^\hbar_{\lambda,s-}H_{B_n}(X^\hbar_{s-},z)\tilde N(\dif s,\dif z).
$$
Thus,  we obtain that for any $t\in[\lambda,T)$,
$$
uM^\hbar_{B_n}(t)+uR^\hbar_{B_n}(t)=\lim_{m\to\infty}\int_\lambda^t\!\!\!\int_{\eps_m<|z|\leq\delta}
uK^\hbar_{\lambda,s-}H_{B_n}(X^\hbar_{s-},z) N(\dif s,\dif z)\stackrel{(\ref{Eq33})}{=}0,
$$
which together with (\ref{Eq4}) yields
$$
\int^t_\lambda uK^\hbar_{\lambda,s}B_{n+1}(X^\hbar_s)\dif s=0,\ \ \forall t\in [\lambda,T).
$$
By the right continuities of $s\mapsto K^\hbar_{\lambda,s}, X^\hbar_s$, we obtain (\ref{Eq5}). The proof is complete.
\end{proof}

\section{Komatsu-Takeuchi's estimate for discontinuous semimartingales}

In this section we shall prove a version of Komatsu-Takeuchi's estimate (cf. \cite{Ko-Ta0, Ko-Ta1}). 
Let $W_t$ be a $d$-dimensional Brownian motion and $N(\dif s,\dif z)$ a Poisson random measure over $\mR^d_0$ 
with intensity measure $\dif s\nu(\dif z)$, where $\nu$ is an infinite measure on $\mR^d_0$ with
$$
\int_{\mR^d_0}(1\wedge|z|^2)\nu(\dif z)<+\infty. 
$$
Let $\tilde N(\dif s,\dif z):=N(\dif s,\dif z)-\dif s\nu(\dif z)$ be the compensated Poisson random measure.

Let $\sS_m$ be the class of $m$-dimensional semimartingales with the following form:
\begin{align}
X_t=\xi_0+\int^t_0\xi^0_s\dif s+\int^t_0\xi^k_s\dif W^k_s+\int^t_0\!\!\!\int_{\mR^d_0}\eta_s(z)\tilde N(\dif s,\dif z),\label{EW1}
\end{align}
where $\xi^k_s,k=0,1,\cdots, d$ and $\eta_s(z)$ are $m$-dimensional predictable processes with
$$
\|X_\cdot(\omega)\|_{\sS_m}:=\sup_{s\in[0,1]}\left(|X_s(\omega)|^2\vee|\xi^0_s(\omega)|^2\vee|\xi^k_s(\omega)|^2
\vee\sup_{z\in\mR^d}\frac{|\eta_s(z,\omega)|^2}{1\wedge|z|^2}\right)<\infty\ \ a.e.-\omega.
$$
Here and below, we use the following convention: If an index appears twice in a product, then it will be summed automatically. For example,
$$
\int^t_0\xi^k_s\dif W^k_s:=\sum_{k=1}^d\int^t_0\xi^k_s\dif W^k_s,\ \ \ |\xi^k_s|^2:=\sum_{k=1}^d|\xi^k_s|^2.
$$
For $\kappa\geq 1$, let $\sS^\kappa_m$ be the subclass of $\sS_m$ with
$$
\|X_\cdot(\omega)\|_{\sS_m}\leq\kappa\ \ a.s.-\omega.
$$

We first prepare the following easy result about the exponential supermartingales.
\bl
Let $(X_t)_{t\geq 0}$ be a one-dimensional semimartingale in $\sS_1$ taking form (\ref{EW1}) with $\xi_0=\xi^0_s=0$. If we let
\begin{align*}
\cE_t(X)&:=\exp\Bigg\{\int^t_0\xi^k_s\dif W^k_s+\int^t_0\!\!\!\int_{\mR^d_0}\eta_s(z)\tilde N(\dif s,\dif z)\\
&-\frac{1}{2}\int^t_0|\xi^k_s|^2\dif s-\int^t_0\!\!\!\int_{\mR^d_0}(\e^{\eta_s(z)}-1-\eta_s(z))\nu(\dif z)\dif s\Bigg\},
\end{align*}
then for any $t\geq 0$,
\begin{align}
\mE \cE_t(X)\leq 1.\label{Eso0}
\end{align}
Moreover, for any $R>0$,  we have on $|R\eta_s(z)|\leq 1$,
\begin{align}
&\left|\int^t_0\xi^k_s\dif W^k_s+\int^t_0\!\!\!\int_{\mR^d_0}\eta_s(z)\tilde N(\dif s,\dif z)-\frac{1}{R}\log\cE_t(RX)\right|\no\\
&\qquad\leq\frac{R}{2}\int^t_0|\xi^k_s|^2\dif s+2R\int^t_0\!\!\!\int_{\mR^d_0}|\eta_s(z)|^2\nu(\dif z)\dif s.\label{Esy}
\end{align}
\el
\begin{proof}
By It\^o's formula, we have
\begin{align*}
\cE_t(X)=1+\int^t_0\cE_s(X)\xi^k_s\dif W^k_s+\int^t_0\!\!\!\int_{\mR^d_0}\cE_{s-}(X)(\e^{\eta_s(z)}-1)\tilde N(\dif s,\dif z).
\end{align*}
Hence,
$$
t\mapsto \cE_t(X)\mbox{ is a positive local martingale (and also a supermartingale)},
$$
which then implies (\ref{Eso0}). As for (\ref{Esy}), it follows by
$$
|\e^x-x-1|\leq 2x^2,  \ \forall |x|\leq 1,
$$
and a direct calculation.
\end{proof}
We are now in a position to prove the following Komatsu-Takeuchi's estimate.
\bt\label{Th42}
For $\kappa\geq 1$, let $(f_t)_{t\geq 0}$ and $(f^0_t)_{t\geq 0}$ be two $m$-dimensional semimartingales in $\sS^\kappa_m$ given by
\begin{align*}
f_t&=f_0+\int^{t\wedge\tau}_0f^0_s\dif s+\int^t_0f^k_s\dif W^k_s+\int^t_0\!\!\!\int_{|z|\leq\delta}g_s(z)\tilde N(\dif s,\dif z),\\
f^0_t&=f^0_0+\int^t_0f^{00}_s\dif s+\int^t_0f^{0k}_s\dif W^k_s+\int^t_0\!\!\!\int_{|z|\leq\delta}g^0_s(z)\tilde N(\dif s,\dif z),
\end{align*}
where $\delta\in(0,1]$ and $\tau$ is a stopping time. For any $\eps, T\in(0,1]$,
there exist positive random variables $\zeta_1$ and $\zeta_2$ with $\mE \zeta_1\leq 1$,
$\mE \zeta_2\leq 1$ such that
\begin{align}
c_1\int^T_0\left\{|f^k_t|^2+\int_{|z|\leq\delta}|g_t(z)|^2\nu(\dif z)\right\}\dif t
\leq(\delta^{-1}+\eps^{-1})\int^T_0|f_t|^2\dif t
+\kappa\delta\log \zeta_1+\kappa(\eps+T\delta)\label{Eso3}
\end{align}
and
\begin{align}
c_2\int^{T\wedge\tau}_0|f^0_t|^2\dif t\leq(\delta^{-\frac{3}{2}}+\eps^{-\frac{3}{2}})\int^T_0|f_t|^2\dif t
+\kappa\delta^{\frac{1}{2}}\log \zeta_2+\kappa(\eps\delta^{-\frac{1}{2}}+\eps^{\frac{1}{2}}+T\delta^{\frac{1}{2}}),\label{Eso2}
\end{align}
where $c_1, c_2\in(0,1)$ only depends on $\int_{|z|\leq 1}|z|^2\nu(\dif z)$.
\et
\begin{proof}
First of all, if $\eps\geq T$, then in view of $f_t\in\sS^\kappa_m$,
$$
\int^T_0\left\{|f^k_s|^2+\int_{|z|\leq\delta}|g_s(z)|^2\nu(\dif z)\right\}\dif t\leq \kappa T\left(1+\int_{|z|\leq 1}|z|^2\nu(\dif z)\right)\leq C\kappa \eps,
$$
and
$$
\int^T_0|f^0_t|^2\dif t\leq \kappa T\leq\kappa\eps,
$$
and hence, (\ref{Eso3}) and (\ref{Eso2}) hold with $\zeta_1=\zeta_2\equiv 1$. Below, 
by replacing $(f_t, f^k_t, g_t(z))$ and $(f^{0k}_t, g^0_t(z))$ with $(f_t, f^k_t, g_t(z))/\sqrt{\kappa}$ and $(f^{0k}_t, g^0_t(z))/\sqrt{\kappa}$,
we shall assume 
$$
\kappa=1,\ \ \eps\in(0,T).
$$
In particular, we have
\begin{align}
|f_t|^2\vee|f^0_t|^2\vee|f^k_t|^2\vee|f^{00}_t|^2\vee|f^{0k}_t|^2
\vee\sup_{z\in\mR^d}\frac{|g_t(z)|^2}{1\wedge|z|^2}\vee\sup_{z\in\mR^d}\frac{|g^0_t(z)|^2}{1\wedge|z|^2}\leq 1.\label{EK2}
\end{align}
We divide the proof into two steps.
\\
\\
{\bf (1)} For $t\in[0,T-\eps]$, by It\^o's formula, we have
\begin{align*}
|f_{t+\eps}|^2&=|f_t|^2+2\int^{t+\eps}_t\<f_s,1_{s\leq\tau}f^0_s\>\dif s+2\int^{t+\eps}_t\<f_s,f^k_s\>\dif W^k_s+\int^{t+\eps}_t|f^k_s|^2\dif s\\
&+\int^{t+\eps}_t\!\!\!\int_{|z|\leq\delta}\Big\{|f_{s-}+g_s(z)|^2-|f_{s-}|^2\Big\}\tilde N(\dif s,\dif z)+\int^{t+\eps}_t\!\!\!\int_{|z|\leq\delta}|g_s(z)|^2\nu(\dif z)\dif s.
\end{align*}
Integrating both sides from $0$ to $T-\eps$ and using Fubini's theorem, we obtain
\begin{align}
\int^{T-\eps}_0|f_{t+\eps}|^2\dif t&\geq 2\int^T_0\eps_s\<f_s,1_{s\leq\tau}f^0_s\>\dif s+\int^T_0\eps_s\left\{|f^k_s|^2+\int_{|z|\leq\delta} |g_s(z)|^2\nu(\dif z)\right\}\dif s\no\\
&+\left\{\int^T_0\xi^{\eps,k}_s\dif W^k_s+\int^T_0\!\!\!\int_{|z|\leq\delta}\eta^\eps_s(z)\tilde N(\dif s,\dif z)\right\}
=:I^\eps_1+I^\eps_2+I^\eps_3,\label{Eso1}
\end{align}
where 
$$
\eps_s:=(T-\eps)\wedge s-0\vee(s-\eps)
$$ 
and
$$
\xi^{\eps,k}_s:=2\eps_s\<f_s,f^k_s\>,\ \ \eta^\eps_s(z):=\eps_s\Big\{|f_{s-}+g_s(z)|^2-|f_{s-}|^2\Big\}.
$$
For $I^\eps_1$,  we have
\begin{align}
|I^\eps_1|\leq 2\eps\int^T_0|f_s||f^0_s|\dif s\leq \int^T_0|f_s|^2\dif s+\eps^2\int^T_0|f^0_s|^2\dif s\stackrel{(\ref{EK2})}{\leq} \int^T_0|f_s|^2\dif s+T\eps^2.\label{EW2}
\end{align}
For $I^\eps_2$, noticing that
$$
|\eps_s-\eps|\leq\eps \Big\{1_{(0,\eps)}(s)+1_{(T-\eps,T)}(s)\Big\},\ \ s\in[0,T],
$$
one easily sees that
\begin{align}
I^\eps_2\stackrel{(\ref{EK2})}{\geq} \eps\int^T_0\left\{|f^k_s|^2+\int_{|z|\leq\delta}|g_s(z)|^2\nu(\dif z)\right\}\dif s-2\left(1+\int_{|z|\leq 1}|z|^2\nu(\dif z)\right)\eps^2.\label{EW3}
\end{align}
For $I^\eps_3$, notice that
$$
|\eta^\eps_s(z)|1_{|z|\leq\delta}\leq \eps(2|f_{s-}||g_s(z)|+|g_s(z)|^2)1_{|z|\leq\delta}\stackrel{(\ref{EK2})}{\leq}3\eps\delta.
$$
Thus, by (\ref{Esy}) with $R=\frac{1}{3\eps\delta}$,
there exists a positive random variable $M^{\eps,\delta}_1$ with $\mE M^{\eps,\delta}_1\leq 1$ such that
\begin{align}
-I^\eps_3&\leq 3\eps\delta\log M^{\eps,\delta}_1+C(\eps\delta)^{-1}\left(\int^T_0|\xi^{\eps,k}_s|^2\dif s
+\int^{T}_0\!\!\!\int_{|z|\leq\delta}|\eta^\eps_s(z)|^2\nu(\dif z)\dif s\right)\no\\
&\!\!\stackrel{(\ref{EK2})}{\leq}3\eps\delta\log M^{\eps,\delta}_1+C\eps\delta^{-1}\left(\int^T_0|f_s|^2\dif s
+\int^{T}_0\!\!\!\int_{|z|\leq\delta}(|f_s|^2|z|^2+|z|^4)\nu(\dif z)\dif s\right)\no\\
&\leq3\eps\delta\log M^{\eps,\delta}_1+C\eps\delta^{-1}\int^T_0|f_s|^2\dif s+CT\eps\delta.\label{Ew4}
\end{align}
Here and below, the constant $C\geq 1$ only depends on $\int_{|z|\leq 1}|z|^2\nu(\dif z)$.
Combining (\ref{Eso1})-(\ref{Ew4}), we obtain
\begin{align}
&\int^T_0\left\{|f^k_s|^2+\int_{|z|\leq\delta}|g_s(z)|^2\nu(\dif z)\right\}\dif s\no\\
&\qquad\leq3\delta\log M^{\eps,\delta}_1+C_1(\delta^{-1}+\eps^{-1})\int^T_0|f_s|^2\dif s+C_1(\eps+T\delta),\label{EW9}
\end{align}
which gives (\ref{Eso3}) with $\kappa=1$ by setting $c_1:=\frac{1}{3\vee C_1}$ and $\zeta_1:=(M^{\eps,\delta}_1)^\frac{3}{3\vee C_1}$.
\\
\\
{\bf (2)} As above, by It\^o's product formula, we have
\begin{align*}
\<f_{t+\eps}, f^0_{t+\eps}\>&=\<f_t, f^0_t\>+\int^{t+\eps}_t\<f_s,f^{00}_s\>\dif s+\int^{t+\eps}_t\<f_s,f^{0k}_s\>\dif W^k_s\\
&\quad+\int^{t+\eps}_t\!\!\!\int_{ |z|\leq\delta}\<f_{s-},g^0_s(z)\>\tilde N(\dif s,\dif z)+\int^{t+\eps}_t1_{s\leq\tau}|f^0_s|^2\dif s\\
&\quad+\int^{t+\eps}_t\<f^k_s,f^0_s\>\dif W^k_s+\int^{t+\eps}_t\!\!\!\int_{ |z|\leq\delta}\<g_s(z),f^0_{s-}\>\tilde N(\dif s,\dif z)\\
&\quad+\int^{t+\eps}_t\<f^k_s,f^{0k}_s\>\dif s+\int^{t+\eps}_t\!\!\!\int_{ |z|\leq\delta}\<g_s(z),g^0_s(z)\>N(\dif s,\dif z).
\end{align*}
Integrating both sides from $0$ to $T-\eps$ and using Fubini's theorem again, we get
\begin{align}
&\int^{T-\eps}_0\{\<f_{t+\eps}, f^0_{t+\eps}\>-\<f_t, f^0_t\>\}\dif t
=\Bigg\{\int^T_0\hat\xi_s^{\eps,k}\dif W^k_s+\int^T_0\!\!\!\int_{ |z|\leq\delta}\hat\eta^\eps_s(z)\tilde N(\dif s,\dif z)\Bigg\}\no\\
&\qquad+\int^T_0\eps_s\left\{\<f_s,f^{00}_s\>+\<f^k_s,f^{0k}_s\>+\int_{|z|\leq\delta}\<g_s(z),g^0_s(z)\>\nu(\dif z)\right\}\dif s\no\\
&\qquad\qquad+\int^T_0\eps_s 1_{s\leq\tau}|f^0_s|^2\dif s=:J^\eps_1+J^\eps_2+J^\eps_3,\label{EW5}
\end{align}
where
$$
\hat\xi^{\eps,k}_s:=\eps_s\Big\{\<f_s,f^{0k}_s\>+\<f^k_s,f^0_s\>\Big\},\ \ \hat\eta^\eps_s(z):=\eps_s\Big\{\<f_{s-},g^0_s(z)\>+\<g_s(z),f^0_{s-}+g^0_s(z)\>\Big\}.
$$
For $J^\eps_1$,  notice that
$$
|\hat\eta^\eps_s(z)|1_{|z|\leq\delta} \leq\eps\Big\{|f_{s-}||g^0_s(z)|+|g_s(z)|(|f^0_{s-}|+|g^0_s(z)|)\Big\}1_{|z|\leq\delta}
\stackrel{(\ref{EK2})}{\leq} 3\eps\delta\leq3\eps\delta^{\frac{1}{2}}.
$$
Thus, by (\ref{Esy}) with $R=\frac{1}{3\eps\delta^{1/2}}$, there exists a positive random variable $M^{\eps,\delta}_2$ with $\mE M^{\eps,\delta}_2\leq 1$ such that
\begin{align}
-J^\eps_1&\leq 3\eps\delta^{\frac{1}{2}}\log M^{\eps,\delta}_2+C(\eps\delta^{\frac{1}{2}})^{-1}\left\{\int^T_0|\hat\xi^{\eps,k}_s|^2\dif s
+\int^{T}_0\!\!\!\int_{|z|\leq\delta}|\hat\eta^\eps_s(z)|^2\nu(\dif z)\dif s\right\}\no\\
%&\leq C\eps\delta\log M^{\eps,\delta}_2+C\eps\delta^{-1}\Bigg\{\int^T_0(|f_s|^2+|f^k_s|^2)\dif s
%+\int^{T}_0\!\!\!\int_{|z|\leq\delta}|f_s|^2|z|^2\nu(\dif z)\dif s+\int^{T}_0\!\!\!\int_{|z|\leq\delta}|g_s(z)|^2\nu(\dif z)\dif s\Bigg\}\no\\
&\!\!\stackrel{(\ref{EK2})}{\leq}3\eps\delta^{\frac{1}{2}}\log M^{\eps,\delta}_2+C\eps\delta^{-\frac{1}{2}}
\int^T_0\left\{|f_s|^2+|f^k_s|^2+\int_{|z|\leq\delta}|g_s(z)|^2\nu(\dif z)\right\}\dif s.\label{EW6}
\end{align}
For $J^\eps_2$, by Young's inequality and (\ref{EK2}), we have
\begin{align*}
|J^\eps_2|&\leq\eps\int^{T}_0\Big\{|f_s||f^{00}_s|+|f^k_s||f^{0k}_s|\Big\}\dif s+\eps\int^T_0\!\!\!\int_{|z|\leq\delta}|g_s(z)||g^0_s(z)|\nu(\dif z)\dif s\no\\
&\leq \eps\delta^{-\frac{1}{2}}\int^T_0\left\{|f_s|^2+|f^k_s|
+\int_{|z|\leq\delta}|g_s(z)|^2\nu(\dif z)\right\}\dif s+CT\eps\delta^{\frac{1}{2}}
\end{align*}
as well as
\begin{align}
&\left|\int^{T-\eps}_0\{\<f_{t+\eps}, f^0_{t+\eps}\>-\<f_t, f^0_t\>\}\dif t\right| \leq\left(\int^T_{T-\eps}+\int^\eps_0\right)|f_tf^0_t|\dif t\no\\
&\quad\leq\eps^{-\frac{1}{2}}\int^T_0|f_t|^2\dif t+C\eps^{\frac{1}{2}}\left(\int^T_{T-\eps}+\int^\eps_0\right)|f^0_t|^2\dif t
\leq\eps^{-\frac{1}{2}}\int^T_0|f_t|^2\dif t+C\eps^{\frac{3}{2}}.\label{EW10}
\end{align}
For $J^\eps_3$,  it is similar to (\ref{EW3}) that
\begin{align}
J^\eps_3\geq\eps\int^T_01_{s\leq\tau}|f^0_s|^2\dif s-2\eps^2=\eps\int^{T\wedge\tau}_0|f^0_s|^2\dif s-2\eps^2.\label{EW7}
\end{align}
Combining (\ref{EW5})-(\ref{EW10}), we obtain
\begin{align*}
\int^{T\wedge\tau}_0|f^0_s|^2\dif s&\leq 3\delta^{\frac{1}{2}}\log M^{\eps,\delta}_2+C(\delta^{-\frac{1}{2}}+\eps^{-\frac{3}{2}})\int^T_0|f_s|^2\dif s
+2\eps+CT\delta^{\frac{1}{2}}+C\eps^{\frac{1}{2}}\no\\
&\quad+C_0\delta^{-\frac{1}{2}}\int^T_0\left\{|f^k_s|^2+\int_{|z|\leq\delta}|g_s(z)|^2\nu(\dif z)\right\}\dif s\no\\
&\stackrel{(\ref{EW9})}{\leq} 3\delta^{\frac{1}{2}}(\log M^{\eps,\delta}_1+C_0\log M^{\eps,\delta}_2)
+C(\delta^{-\frac{3}{2}}+\eps^{-1}\delta^{-\frac{1}{2}}+\eps^{-\frac{3}{2}})\int^T_0|f_s|^2\dif s\no\\
&\quad+2\eps+CT\delta^{\frac{1}{2}}+C\eps^{\frac{1}{2}}+C(T\delta+\eps)\delta^{-\frac{1}{2}},
\end{align*}
which yields (\ref{Eso2}) with $\kappa=1$ by suitable choices of $c_2$ and $\zeta_2$ as in {\bf (1)}. The proof is thus complete.
\end{proof}

%\bc In the situation of Theorem \ref{Th42}, if the L\'evy measure $\nu$ satisfies the following order condition: for some $\alpha\in(0,2)$ and positive definite matrix $\Lambda$,
%$$ \inf_{\delta\in(0,1)}\frac{1}{\delta^{2-\alpha}}\int_{|z|\leq\delta}|\<z,x\>|^2\nu(\dif z)\geq |\Lambda x|^2,$$
%and $z\mapsto g_s(z)\in C^2(\mR^d)$ satisfies $$|\nabla^2_z g_s(z)|\leq\kappa,$$then for any $\delta, T\in(0,1]$, there exist positive random variables $\zeta$ with $\mE \zeta\leq 1$ such that\begin{align}c_0\int^T_0\left\{|f^k_t|^2+|\Lambda^{k\cdot}\cdot \nabla g_t(0)|^2\right\}\dif t\leq(\delta^{-1}+\eps^{-1})\int^T_0|f_t|^2\dif t+\kappa\delta\log M^{\eps,\delta}_T+\kappa(\eps+T\delta)\label{Eso3}\end{align}and\begin{align}c_0\int^T_0|f^0_t|^2\dif t\leq(\delta^{-\frac{3}{2}}+\eps^{-\frac{3}{2}})\int^T_0|f_t|^2\dif t+\kappa\delta^{\frac{1}{2}}\log \tilde M^{\eps,\delta}_T+\kappa(\eps\delta^{-\frac{1}{2}}+\eps^{\frac{1}{2}}+T\delta^{\frac{1}{2}}),\label{Eso2}\end{align}where $c_0\in(0,1)$ only depends on $\int_{|z|\leq 1}|z|^2\nu(\dif z)$.\ec\begin{proof}Noticing that$$g_t(z)=z\cdot\nabla g_t(0)+\gamma_t(z)$$with $|\gamma_t(z)|^2\leq \kappa|z|^4$, by $|a+b|^2\geq\frac{|a|^2}{2}-|b|^2$, we have\begin{align*}\int_{|z|\leq\delta}|g_t(z)|^2\nu(\dif z)&\geq\frac{1}{2}\int_{|z|\leq\delta}|z\cdot \nabla g_t(0)|^2\nu(\dif z)-\int_{|z|\leq\delta}|\gamma_t(z)|^2\nu(\dif z)\\&\geq\frac{\delta^{2-\alpha}}{2}|\Lambda^{k\cdot}\cdot \nabla g_t(0)|^2-\kappa\int_{|z|\leq\delta}|z^4\nu(\dif z)\end{align*}\end{proof}
\br
The factor $T$ before $\delta^{\frac{1}{2}}$ in (\ref{Eso2}) is quite important for the proof of Lemma \ref{Le9} below.
Meanwhile, the freedom of $\eps$ allows it to be any positive power of $\delta$ so that for any $m\geq 1$ and $T\in(0,1)$,
\begin{align}
\frac{c_0}{2\kappa}\int^T_0|f^0_t|^2\dif t\leq\delta^{-\frac{3m}{2}}\int^T_0|f_t|^2\dif t
+\delta^{\frac{1}{2}}\log \zeta^{\frac{1}{2}}_2+(\delta^{\frac{m}{2}}+T\delta^{\frac{1}{2}}).\label{Eso22}
\end{align}
Moreover, the role of stopping time $\tau$ will become clear in the following corollary.
\er
We have the following corollary.
\bc
Let $(f_t)_{t\geq 0}$ and $(f^0_t)_{t\geq 0}$ be two $m$-dimensional semimartingales in $\sS_m$ as given in Theorem \ref{Th42}.
For $\delta\in(0,1)$ and $\gamma\in(0,\frac{1}{4})$, define a stopping time:
$$
\tau:=\inf\left\{t\geq 0: |f_t|^2\vee|f^0_t|^2\vee|f^k_t|^2\vee|f^{00}_t|^2\vee|f^{0k}_t|^2
\vee\sup_{z\in\mR^d}\frac{|g_t(z)|^2}{1\wedge|z|^2}\vee\sup_{z\in\mR^d}\frac{|g^0_t(z)|^2}{1\wedge|z|^2}\geq\delta^{-\gamma}\right\}.
$$
For any $\delta, T\in(0,1]$ and $\gamma\in(0,\frac{1}{4})$, we have
$$
\mP\left\{\int^T_0|f_t|^2\dif t\leq \delta^{2-\gamma},\int^T_0|f^0_t|^2\dif t\geq \delta^{\frac{1}{4}-\gamma}; \tau\geq T\right\}
\leq \e^{5-c_0\delta^{-\frac{1}{4}}}.
$$
\ec
\begin{proof}
In Theorem \ref{Th42}, if we replace $f_t, f^k_t, g_t(z)$ and $f^0_t,  f^{0k}_t, g^0_t(z)$ by $f_{t\wedge\tau}, 1_{t\leq\tau}f^k_t, 1_{t\leq\tau}g_t(z)$ 
and $f^0_{t\wedge\tau},  1_{t\leq\tau}f^{0k}_t, 1_{t\leq\tau}g^0_t(z)$ respectively, then
by (\ref{Eso2}) with $\eps=\delta$ and $\kappa=\delta^{-\gamma}$, we obtain
$$
c_0\int^{T\wedge\tau}_0|f^0_t|^2\dif t\leq 2\delta^{-\frac{3}{2}}\int^T_0|f_{t\wedge\tau}|^2\dif t+\delta^{\frac{1}{2}-\gamma}\log \zeta_2+3\delta^{\frac{1}{2}-\gamma}.
$$
In particular,
\begin{align*}
\left\{\int^T_0|f_t|^2\dif t\leq \delta^{2-\gamma},\int^T_0|f^0_t|^2\dif t\geq \delta^{\frac{1}{4}-\gamma}; \tau\geq T\right\}
\subset\Big\{\zeta_2\geq \e^{c_0\delta^{-\frac{1}{4}}-5}\Big\},
\end{align*}
which gives the desired estimate by Chebyshev's inequality and $\mE\zeta_2\leq 1$.
\end{proof}
\br
This type of Norris' estimate was proven by Cass in \cite{Ca}. Our proof is simpler, and the dependence of constants are more explicit.
\er

\section{Proof of  Theorem \ref{Th1}}

In this section we assume that
$$
\sigma(x)=\sigma\mbox{ is constant and $b\in C^\infty_b(\mR^d)$}.
$$
Let us consider the following SDE:
$$
X_t(x)=x+\int^t_0b(X_s(x))\dif s+\sigma W_{S_t}.
$$
Let $J_t(x)=\nabla X_t(x)$ be the Jacobian matrix of $X_t(x)$. It is easy to see that
$$
J_t(x)=\mI+\int^t_0\nabla b(X_s(x)) J_s(x)\dif s.
$$
Let $K_t(x)$ be the inverse matrix of $J_t(x)$. Then
$$
K_t(x)=\mI-\int^t_0K_s(x)\nabla b(X_s(x))\dif s.
$$
In this case, since $J_{s,t}(x)=J_t(x)K_s(x)$, the Malliavin covariance matrix $\Sigma_t(x)$ given by (\ref{Ep4}) can be written as:
\begin{align}
\Sigma_t(x)=J_t(x)\left(\int^t_0K_s(x)\sigma\sigma^*K^*_s(x)\dif S_s\right)J^*_t(x).\label{Ep44}
\end{align}
The core task of this section is to prove the $L^p$-integrability of the inverse of $\Sigma_t(x)$.

Below, for simplicity of notations, we shall drop the variable ``$x$'' unless necessary. We first prove the following key lemma.
\bl\label{Le9}
Let $V:\mR^d\to\mR^d\times\mR^d$ be a bounded smooth function with bounded derivatives of all orders.   Under {\bf(AS$^\alpha$)},
for any $\beta\in(0\vee(4\alpha-7),1)$, there exist two constants $C_1\geq 1$ and $c_1\in(0,1)$ only depending on $b,V$ and $\alpha,\beta,\nu_L$ such that for all $\delta, t\in(0,1)$,
\begin{align}
\sup_{|u|=1}\mP\left(\int^t_0|uK_s[b,V](X_s)|^2\dif s\geq t\delta^{\frac{1-\beta}{2}},\ \int^t_0|uK_sV(X_s)|^2\dif s\leq t\delta^{9-\frac{\beta}{2}}\right)
\leq C_1\exp\{-c_1t\delta^{-\frac{\beta}{2}}\},\label{EE4}
\end{align}
where $[b,V]:=b\cdot\nabla V-V\cdot\nabla b$.
\el
\begin{proof}
We divide the proof into three steps.
\\
{\bf (1)} Let $\hbar$ be a c\`adl\`ag purely discontinuous $\mR^d$-valued function with finitely many jumps and $\hbar_0=0$.
Let $\ell$ be another c\`adl\`ag purely discontinuous $\mR^d$-valued function with $\ell_0=0$.
Let $X_t(x;\ell,\hbar)$ solve the following equation:
$$
X_t(x;\ell,\hbar)=x+\int^t_0b(X_s(x;\ell,\hbar))\dif s+\sigma\ell_t+\sigma\hbar_t.
$$
Let $K_t(x;\ell,\hbar)$ solve the following matrix-valued ODE:
\begin{align}
K_t(x;\ell,\hbar)=\mI-\int^t_0K_s(x;\ell,\hbar)\nabla b(X_s(x;\ell,\hbar))\dif s.\label{EK3}
\end{align}
By Gronwall's inequality, it is clear that
\begin{align}
|K_t(x;\ell,\hbar)|\leq \e^{\|\nabla b\|_\infty},\ \ \forall t\in[0,1].\label{EC1}
\end{align}
For $\delta\in(0,1)$, let $L^\delta_t$ and $\hat L^\delta_t$ be defined as in (\ref{NB4}). Write
$$
X^\delta_t(x):=X_t(x;L^\delta_\cdot,0),\ \  K^\delta_t(x):=K_t(x;L^\delta_\cdot,0).
$$
Then $X^\delta_t(x)$ solves the following SDE:
\begin{align}
X^\delta_t(x)=x+\int^t_0 b(X^\delta_s(x))\dif s+\int_{|z|\leq\delta}\sigma z\tilde N((0,t],\dif z).\label{EK33}
\end{align}
Define
\begin{align*}
V_0(x)&:=[b,V](x)+\int_{|z|\leq\delta}(V(x+\sigma z)-V(x)-\sigma z\cdot \nabla V(x))\nu_L(\dif z),\\
V_1(x)&:=[b,V_0](x)+\int_{|z|\leq\delta}(V_0(x+\sigma z)-V_0(x)-\sigma z\cdot \nabla V_0(x))\nu_L(\dif z).
\end{align*}
For any $u\in\mR^d$ with $|u|\leq \e^{\|\nabla b\|_\infty}$, set
$$
f_t:=uK^\delta_t V(X^\delta_t),\  \ g_t(z):=uK^\delta_t(V(X^\delta_t+\sigma z)-V(X^\delta_t))
$$
and
$$
f^0_t:=uK^\delta_t V_0(X^\delta_t),\  f^{00}_t:=uK^\delta_t V_1(X^\delta_t),\ g^0_t(z):=uK^\delta_t(V_0(X^\delta_t+\sigma z)-V_0(X^\delta_t)).
$$
By equations (\ref{EK3}), (\ref{EK33}) and using It\^o's formula (see (\ref{IT0})), we have
\begin{align*}
f_t&=u V(x)+\int^t_0uK^\delta_s[b, V](X^\delta_s)\dif s+\int^t_0\!\!\!\int_{|z|\leq\delta}g_{s-}(z)\tilde N(\dif s,\dif z)\\
&+\int^t_0\!\!\!\int_{|z|\leq\delta}uK^\delta_s(V(X^\delta_s+\sigma z)-V(X^\delta_s)-\sigma z\cdot \nabla V(X^\delta_s))\nu_L(\dif z)\dif s\\
&=u V(x)+\int^t_0f^0_s\dif s+\int^t_0\!\!\!\int_{|z|\leq\delta}g_{s-}(z)\tilde N(\dif s,\dif z)
\end{align*}
and
$$
f^0_t=uV_0(x)+\int^t_0f^{00}_s\dif s+\int^t_0\!\!\!\int_{|z|\leq\delta}g^0_{s-}(z)\tilde N(\dif s,\dif z).
$$
Since $V$ has bounded derivatives of all orders, by (\ref{EC1}) and $|u|\leq\e^{\|\nabla b\|_\infty}$, there exists a constant $\kappa\geq 1$ only depending on $b$ and $V$ such that
for all $t\in[0,1]$ and $z\in\mR^d$,
$$
|f_t|^2,\ |f^0_t|^2,\ |f^{00}_t|^2\leq \kappa,\ \ |g_t(z)|^2,\ |g^0_t(z)|^2\leq \kappa(1\wedge|z|^2).
$$
By (\ref{Eso22}) with $m=5$, for any $t\in(0,1)$, there exists a positive random variable $\zeta$ with $\mE\zeta\leq 1$ such that
\begin{align*}
2c_0\int^t_0|f^0_s|^2\dif s&\leq \delta^{-\frac{15}{2}}\int^t_0|f_s|^2\dif s+\delta^{\frac{1}{2}}\log \zeta
+(\delta^{\frac{5}{2}}+t\delta^{\frac{1}{2}})\ \ a.s.,
\end{align*}
where $c_0\in(0,1)$  only depends on $\kappa$ and $\int_{|z|\leq 1}|z|^2\nu_L(\dif z)$. From this, dividing both sides by $\delta^{\frac{1}{2}}$ and
taking exponential and expectations, we derive that
\begin{align}
\mE\left(\exp\left\{2c_0\delta^{-\frac{1}{2}}\int^t_0|f^0_s|^2\dif s-\delta^{-8}\int^t_0|f_s|^2\dif s\right\}\right)
\leq\exp\{\delta^2+t\}.\label{EC2}
\end{align}
Recalling the definition of $f^0_t$ and by $|x+y|^2\geq\frac{|x|^2}{2}-|y|^2$, we have
\begin{align*}
|f^0_t|^2&\geq \frac{|uK^\delta_t[b,V](X^\delta_t)|^2}{2}-\e^{4\|\nabla b\|_\infty}\|\nabla^2 V\|^2_\infty|\sigma|^4\left(\int_{|z|\leq\delta}|z|^2\nu_L(\dif z)\right)^2\\
&\stackrel{(\ref{EUT0})}{\geq} \frac{|uK^\delta_t[b,V](X^\delta_t)|^2}{2}-C_2\delta^{4-2\alpha}.
\end{align*}
Thus,  by (\ref{EC2}), we obtain that for all $\delta,t\in(0,1)$,
\begin{align}
&\sup_{|u|\leq\e^{\|\nabla b\|_\infty}}\mE\left(\exp\left\{c_0\delta^{-\frac{1}{2}}\int^t_0|uK^\delta_s[b,V](X^\delta_s)|^2\dif s
-\delta^{-8}\int^t_0|uK^\delta_sV(X^\delta_s)|^2\dif s\right\}\right)\no\\
&\qquad\qquad\qquad\leq\exp\Big\{\delta^2+t(C_3\delta^{\frac{7}{2}-2\alpha}+1)\Big\},\label{EC3}
\end{align}
where $C_3:=2C_2c_0$.
\\
\\
{\bf (2)} Let $\hbar$ be a c\`adl\`ag purely discontinuous $\mR^d$-valued function with finitely many jumps and $\hbar_0=0$.
Write for $\delta\in(0,1)$,
$$
X^\delta_t(x;\hbar):=X_t(x;L^\delta_\cdot,\hbar),\ \ K^\delta_t(x;\hbar):=K_t(x;L^\delta_\cdot,\hbar)
$$
and for $u\in\mR^d$,
\begin{align}
\cJ^\delta_t(u,x;\hbar):=\exp\Bigg\{&c_0\delta^{-\frac{1}{2}}\int^t_0|uK^\delta_s(x;\hbar)[b,V](X^\delta_s(x;\hbar))|^2\dif s\no\\
&-\delta^{-8}\int^t_0|uK^\delta_s(x;\hbar)V(X^\delta_s(x;\hbar))|^2\dif s\Bigg\}.\label{Def0}
\end{align}
By (\ref{EC3}), we have
\begin{align}
\sup_{x\in\mR^d}\sup_{|u|\leq\e^{\|\nabla b\|_\infty}}\mE\cJ^\delta_t(u,x;0)\leq\exp\Big\{\delta^2+t(C_3\delta^{\frac{7}{2}-2\alpha}+1)\Big\}.\label{EC4}
\end{align}
%Clearly, $X^\delta_t(x;\hbar)$ solves the following SDE:
%$$
%X^{\delta}_t(x;\hbar)=x+\int^t_0b(X^\delta_s(x;\hbar))\dif s+\sigma L^\delta_t+\sigma\hbar_t
%$$
%and
%$$
%K^\delta_t(x;\hbar)=(\nabla X^{\delta}_t(x;\hbar))^{-1}.
%$$
For $t\in(0,1)$, let $n$ be the jump number of $\hbar$ before time $t$. Let $0=:t_0<t_1<\cdots<t_n\leq t=:t_{n+1}$ be the jump time of $\hbar$. It is easy to see that
for any $s\in[0,t_{j+1}-t_j)$,
\begin{align}
X^\delta_{s+t_j}(x;\hbar)=X_s(X^\delta_{t_j}(x;\hbar);\vartheta_{t_j}L^\delta_\cdot,0)\label{ETY1}
\end{align}
and
\begin{align}
K^\delta_{s+t_j}(x;\hbar)=K^\delta_{t_j}(x;\hbar)\cdot K_s(X^\delta_{t_j}(x;\hbar);\vartheta_{t_j}L^\delta_\cdot,0),\label{ETY11}
\end{align}
where
$$
(\vartheta_{t_j}L^\delta_s:=L^\delta_{s+t_j}-L^\delta_{t_j})_{s\geq 0}\stackrel{(d)}{=}(L^\delta_s)_{s\geq 0}.
$$
Noticing that for $j=1,\cdots,n$,
%\begin{align*}
%\cJ^{\delta}_{t_{j+1}}(u,x;\hbar)=\cJ^{\delta}_{t_j}(u,x;\hbar)\cdot\cJ^{\delta}_{t_{j+1}-t_j}(uK^\delta_{t_j}(x;\hbar),X^\delta_{t_j}(x;\hbar);0)
%\end{align*}
$$
\mbox{$\cJ^{\delta}_{t_j}(u,x;\hbar)$, $K^\delta_{t_j}(x;\hbar)$, $X^\delta_{t_j}(x;\hbar)$ are independent of $X_\cdot(\cdot;\vartheta_{t_j}L^\delta_\cdot,0)$
and $K_\cdot(\cdot;\vartheta_{t_j}L^\delta_\cdot,0)$}
$$
and
$$
X_\cdot(\cdot;\vartheta_{t_j}L^\delta_\cdot,0)\stackrel{(d)}{=}X^\delta_s(\cdot;0),\ \ K_\cdot(\cdot;\vartheta_{t_j}L^\delta_\cdot,0)\stackrel{(d)}{=}K^\delta_s(\cdot;0),
$$
by (\ref{ETY1}) and (\ref{ETY11}), we obtain that for $|u|=1$,
\begin{align}
\mE \cJ^{\delta}_{t_{n+1}}(u,x;\hbar)&=\mE\left(\cJ^{\delta}_{t_n}(u,x;\hbar)
\cdot\Big(\mE\cJ^{\delta}_{t_{n+1}-t_n}(u',y;0)\Big)\Big|_{u'=uK^\delta_{t_n}(x;\hbar), y=X^\delta_{t_n}(x;\hbar)}\right)\no\\
&\!\!\!\!\!\!\stackrel{(\ref{EC1}) (\ref{EC4})}{\leq}\exp\Big\{\delta^2+(t_{n+1}-t_n)(C_3\delta^{\frac{7}{2}-2\alpha}+1)\Big\}\mE\cJ^{\delta}_{t_n}(u,x;\hbar)\no\\
&\leq\cdots\cdots\cdots\cdots\cdots\cdots\no\\
&\leq\Pi_{j=0}^n\exp\Big\{\delta^2+(t_{j+1}-t_j)(C_3\delta^{\frac{7}{2}-2\alpha}+1)\Big\}\no\\
&=\exp\Big\{\delta^2(n+1)+t_{n+1}(C_3\delta^{\frac{7}{2}-2\alpha}+1)\Big\}.\label{LK1}
\end{align}
Let $N^\delta_t$ be the jump number of $\hat L^\delta_\cdot$ before time $t$, i.e.,
$$
N^\delta_t=\sum_{s\in(0,t]}1_{|\Delta\hat L^\delta_t|>0}=\int_{|z|>\delta}N((0,t],\dif z)=\sum_{s\in(0,t]}1_{|\Delta L_s|>\delta},
$$
which is a Poisson process with intensity $\int_{|z|>\delta}\nu_L(\dif z)=:\lambda_\delta$.
Observing that
$$
X_\cdot(x)=X^\delta_\cdot(x; \hat L^\delta),\ \ K_\cdot(x)=K^\delta_\cdot(x; \hat L^\delta),
$$
and recalling (\ref{Def0}) and the independence of $L^\delta$ and $\hat L^\delta$, we have for any $|u|=1$,
\begin{align*}
&\mE\left( \exp\left\{c_0\delta^{-\frac{1}{2}}\int^t_0|uK_s[b,V](X_s)|^2\dif s-\delta^{-8}\int^t_0|uK_sV(X_s)|^2\dif s\right\}\right)\\
&\quad=\mE\Big(\mE\cJ^\delta_t(u,x;\hbar)|_{\hbar=\hat L^\delta}\Big)
=\sum_{n=0}^\infty\mE\left(\Big(\mE\cJ^{\delta}_t(u,x;\hbar)\Big)_{\hbar=\hat L^\delta}; N^\delta_t=n\right)\\
&\quad\stackrel{(\ref{LK1})}{\leq} \sum_{n=0}^\infty\exp\Big\{\delta^2(n+1)+t(C_3\delta^{\frac{7}{2}-2\alpha}+1)\Big\}\mP(N^\delta_t=n)\\
&\quad=\exp\Big\{\delta^2+t(C_3\delta^{\frac{7}{2}-2\alpha}+1)\Big\}\sum_{n=0}^\infty\e^{\delta^2n}\frac{(t\lambda_\delta)^n}{n!}\e^{-t\lambda_\delta}\\
&\quad=\exp\Big\{\delta^2+t(C_3\delta^{\frac{7}{2}-2\alpha}+1)+(\e^{\delta^2}-1)t\lambda_\delta\Big\}\\
&\quad\stackrel{(\ref{EUT0})}{\leq} \exp\Big\{2+C_4t\delta^{\frac{7}{2}-2\alpha}+C_5t\delta^{2-\alpha}\Big\},\ \ \forall t,\delta\in(0,1),
\end{align*}
where in the last step we have used that $\e^x-1\leq 3x$ for $x\in(0,1)$.
\\
\\
{\bf (3)} By Chebyshev's inequality, we have for any $\beta\in(0,1)$,
\begin{align*}
&\sup_{|u|=1}\mP\left\{c_0\delta^{-\frac{1}{2}}\int^t_0|uK_s[b,V](X_s)|^2\dif s-\delta^{-8}\int^t_0|uK_sV(X_s)|^2\dif s\geq t\delta^{-\frac{\beta}{2}}\right\}\no\\
&\qquad\qquad\leq\exp\Big\{2+C_4t\delta^{\frac{7}{2}-2\alpha}+C_5t\delta^{2-\alpha}-t\delta^{-\frac{\beta}{2}}\Big\}.
\end{align*}
In particular, if $\beta\in(0\vee(4\alpha-7),1)$, then there exists a constant $\delta_0$ such that for all $\delta\in(0,\delta_0)$ and $t\in(0,1)$,
$$
\sup_{|u|=1}\mP\left\{\int^t_0|uK_s[b,V](X_s)|^2\dif s\geq \frac{2t\delta^{\frac{1-\beta}{2}}}{c_0},\int^t_0|uK_sV(X_s)|^2\dif s\leq t\delta^{8-\frac{\beta}{2}}\right\}
\leq\exp\{3-t\delta^{-\frac{\beta}{2}}\},
$$
which then gives the desired estimate by adjusting the constants and rescaling $\delta$.
\end{proof}
\bl
Let $t_0:=1\wedge(\frac{\e^{-\|\nabla b\|_\infty}}{2\|\nabla b\|_\infty})$.
Under (\ref{UH})  and {\bf(AS$^\alpha$)}, there exist $\gamma=\gamma(\alpha)\in(0,1)$, $C_2\geq 1$ and $c_2\in(0,1)$ such that for all $t\in(0,t_0)$ and $\lambda\geq 1$,
\begin{align}
\sup_{|u|=1}\sup_{x\in\mR^d}\mE\exp\left\{-\lambda\int^t_0|uK_s(x)\sigma|^2\dif s\right\}\leq C_2\exp\{- c_2 t\lambda^\gamma\}.\label{EE8}
\end{align}
\el
\begin{proof}
For fixed $\beta\in(0\vee(4\alpha-7),1)$, set $a:=\frac{1-\beta}{18-\beta}$ and define for $j=1,\cdots, j_0$, 
$$
E_j:=\left\{\int^t_0|uK_sB_j(X_s)|^2\dif s\leq t\delta^{a^j(9-\frac{\beta}{2})}\right\}.
$$
Since $a^{j+1}(9-\frac{\beta}{2})=\frac{a^j(1-\beta)}{2}$ and $B_{j+1}=[b,B_j]$, by (\ref{EE4}) with $\delta$ replaced by $\delta^{a^j}$, we have
\begin{align}
\mP(E_j\cap E_{j+1}^c)\leq C_1\exp\{-c_1t\delta^{-a^j\frac{\beta}{2}}\}.\label{EC5}
\end{align}
Noticing that
$$
E_1\subset\left(\cap_{j=1}^{j_0} E_j\right)\cup\left(\cup_{j=1}^{j_0-1}(E_j\cap E_{j+1}^c)\right),
$$
we have
\begin{align}
\mP(E_1)\leq \mP\left(\cap_{j=1}^{j_0} E_j\right)+\sum_{j=1}^{j_0-1}\mP(E_j\cap E_{j+1}^c).\label{EC6}
\end{align}
On the other hand, noticing that for $t\leq 1\wedge(\frac{\e^{-\|\nabla b\|_\infty}}{2\|\nabla b\|_\infty})=:t_0$,
$$
\inf_{|u|=1}|uK_t|\geq 1-|K_t-\mI|\geq 1-\|\nabla b\|_\infty \e^{\|\nabla b\|_\infty} t\geq\tfrac{1}{2},
$$
by (\ref{UH}) we have
\begin{align}
\bigcap_{j=1}^{j_0}E_j&\subset\left\{\sum_{j=1}^{j_0}\int^t_0|uK_sB_j(X_s)|^2\dif s\leq t\sum_{j=1}^{j_0}\delta^{a^j(9-\frac{\beta}{2})}\right\}\no\\
&\subset\left\{\kappa_1\int^t_0|uK_s|^2\dif s\leq t \sum_{j=1}^{j_0}\delta^{a^{j_0}(9-\frac{\beta}{2})}\right\}\no\\
&\subset\left\{\frac{t\kappa_1}{2}\leq t j_0\delta^{a^{j_0}(9-\frac{\beta}{2})}\right\}=\emptyset,\label{EC7}
\end{align}
provided $t\leq t_0$ and $\delta\leq\delta_1$ for some $\delta_1$ small enough.
Therefore, combining (\ref{EC5})-(\ref{EC7}), we obtain that for all $t\in(0,t_0)$ and $\delta\in(0,\delta_1)$,
$$
\mP\left\{\int^t_0|uK_s\sigma|^2\dif s\leq t\delta^{a(9-\frac{\beta}{2})}\right\}\leq C_1 j_0\exp\Big\{-c_1t\delta^{-a^{j_0}\frac{\beta}{2}}\Big\}.
$$
which in turn  implies that for $\theta=\frac{\beta a^{j_0}}{1-\beta}$ and all $\eps\in(0,1)$, $t\in(0,t_0)$,
$$
\mP\left\{\int^t_0|uK_s\sigma|^2\dif s\leq t\eps\right\}\leq C\exp\Big\{-c_1t\eps^{-\theta}\Big\}.
$$
For $\lambda\geq t$, setting $r:=(\lambda/t)^{\frac{-1}{1+\theta}}$ and $\xi:=\frac{1}{t}\int^t_0|uK_s\sigma|^2\dif s$, we have
\begin{align*}
\mE\e^{-\lambda\xi}&=\int^\infty_0\lambda\e^{-\lambda\eps}\mP(\xi\leq\eps)\dif \eps\\
&\leq\int^\infty_r\lambda\e^{-\lambda\eps}\dif\eps+C\int^r_0\lambda\e^{-\lambda\eps-c_1t\eps^{-\theta}}\dif \eps\\
&=\e^{-\lambda r}+C\int^{\lambda r}_0\e^{-s-c_1t \lambda^\theta s^{-\theta}}\dif s\\
&\leq\e^{-\lambda r}+C\e^{-c_1t r^{-\theta}}\int^{\lambda r}_0\e^{-s}\dif s\\
&\leq\e^{-t (\lambda/t)^{\frac{\theta}{1+\theta}}}+C\e^{-c_1t (\lambda/t)^{\frac{\theta}{1+\theta}}}.
\end{align*}
By replacing $\lambda$ with $\lambda t$, we obtain the desired estimate with $\gamma=\frac{\theta}{1+\theta}$.
\end{proof}

\bl\label{Le4}
Let $t_0:=1\wedge(\frac{\e^{-\|\nabla b\|_\infty}}{2\|\nabla b\|_\infty})$. Under (\ref{UH})  and {\bf(AS$^\alpha$)}, if $S_t$ has finite moments of all orders, 
then for any $p\geq 1$,  there exists a constant $C\geq 1$ such that for all $t\in(0,t_0)$,
\begin{align}
\sup_{x\in\mR^d}\mE \Big((\det\Sigma_t(x))^{-p}\Big)\leq Ct^{-\frac{2(p+1)}{\alpha\gamma}},\label{EK9}
\end{align}
where $\gamma$ is the same as in (\ref{EE8}).
\el
\begin{proof}
Write
$$
\tilde\Sigma_t:=\int^t_0K_s\sigma \sigma^*K^*_s\dif S_s.
$$
By Lemma \ref{Le10} and (\ref{EE8}), we have for all $t\in(0,t_0)$ and $\lambda\geq 1$,
\begin{align}
\mE\exp\{-\lambda (u\tilde\Sigma_tu^*)\}&=\mE\exp\left\{-\lambda\int^t_0|uK_s\sigma|^2\dif S_s\right\}\no\\
&\leq\left(\mE\exp\left\{-c_0\lambda^{\frac{\alpha}{2}} \int^t_0|uK_s\sigma|^2\dif s\right\}\right)^{\frac{1}{2}}\no\\
&\leq C^{\frac{1}{2}}_2\exp\Big\{-c_2t(c_0\lambda^{\frac{\alpha}{2}})^\gamma/2\Big\}.\label{EK7}
\end{align}
Set
$$
F(u;\lambda):=\exp\{-\lambda (u\tilde\Sigma_tu^*)\}.
$$
For any $u_0,u_1\in\mS^{d-1}$, let $(u_r)_{r\in(0,1)}$ be the shortest arc connecting $u_0$ and $u_1$. By the mean value theorem, one sees that for some $r_0\in(0,1)$,
$$
|F(u_0;\lambda)-F(u_1;\lambda)|\leq 2\rho(u_0,u_1) F(u_{r_0};\lambda)\lambda|\tilde\Sigma_t|,
$$
where $\rho(u_0,u_1)$ denotes the distance of sphere. Since
$$
|\tilde\Sigma_t|\leq C|S_t|,
$$
by H\"older's inequality and (\ref{EK7}), there exist constants $C_3\geq 1, c_3\in(0,1)$
such that for all $t\in(0,t_0)$, $\lambda\geq 1$ and $u_0,u_1\in\mS^{d-1}$,
$$
\mE|F(u_0;\lambda)-F(u_1;\lambda)|^d\leq C_3\rho(u_0,u_1)^d\exp\{-c_3t\lambda^{\frac{\alpha\gamma}{2}}\}\lambda^d.
$$
Since $\mS^{d-1}$ is compact and has dimension $d-1$, by Kolmogorov's continuity theorem, we have
$$
\mE\left(\sup_{|u|=1}F(u;\lambda)\right)\leq C_4\exp\{-c_4t\lambda^{\frac{\alpha\gamma}{2}}\}\lambda,\ \ t\in(0,t_0),\ \lambda\geq 1.
$$
Noticing that for a $d\times d$-symmetric matrix $A$, 
$$
(\det A)^{\frac{1}{d}}\geq\inf_{|u|=1}uAu^*,
$$
we have
\begin{align*}
\mE \Big((\det\tilde\Sigma_t)^{-\frac{p}{d}}\Big)&=\frac{1}{\Gamma(p)}\int^\infty_0\lambda^{p-1}\mE\exp\Big\{-\lambda(\det\tilde\Sigma_t)^{\frac{1}{d}}\Big\}\dif\lambda\\
&\leq\frac{1}{\Gamma(p)}\int^\infty_0\lambda^{p-1}\mE\exp\Big\{-\lambda\inf_{|u|=1} (u\tilde\Sigma_t u^*)\Big\}\dif\lambda\\
&=\frac{1}{\Gamma(p)}\int^\infty_0\lambda^{p-1}\mE\left(\sup_{|u|=1}F(u;\lambda)\right)\dif\lambda\\
&\leq\frac{1}{\Gamma(p)}\left(\int^1_0\lambda^{p-1}\dif\lambda+C_4\int^\infty_1\lambda^p
\exp\{-c_4t\lambda^{\frac{\alpha\gamma}{2}}\}\dif\lambda\right)\\
&=\frac{1}{\Gamma(p)}\left(\frac{1}{p}+C_4t^{-\frac{2(p+1)}{\alpha\gamma}}
\int^\infty_{t^{2/(\alpha\gamma)}}\lambda^p\exp\{-c_4\lambda^{\frac{\alpha\gamma}{2}}\}\dif\lambda\right)\\
&\leq\frac{1}{\Gamma(p)}\left(\frac{1}{p}+C_5t^{-\frac{2(p+1)}{\alpha\gamma}}\right),
\end{align*}
which gives 
\begin{align}
\mE \Big((\det\tilde\Sigma_t)^{-\frac{p}{d}}\Big)\leq C_6t^{-\frac{2(p+1)}{\alpha\gamma}},\ \ t\in(0,t_0).\label{EK8}
\end{align}
On the other hand, by (\ref{Ep44}) we have
$$
\Sigma_t=J_t\tilde\Sigma_t J^*_t.
$$
Since
$$
\det J_t=\exp\left\{\int^t_0(\div b)(X_s)\dif s\right\}\geq\exp\{-t\|\div b\|_\infty\},
$$
which together with (\ref{EK8}) yields (\ref{EK9}).
\end{proof}

\bt\label{Th4}
Let $t_0:=1\wedge(\frac{\e^{-\|\nabla b\|_\infty}}{2\|\nabla b\|_\infty})$. 
Under (\ref{UH}) and {\bf(AS$^\alpha$)}, for any $k,m\in\mN_0$ with $k+m\geq 1$, there are $\gamma_{k,m}>0$
and $C=C(k,m)>0$ such that for all $f\in C^\infty_b(\mR^d)$ and $t\in(0,t_0)$,
\begin{align}
\sup_{x\in\mR^d}\Big|\nabla^k\mE\Big( (\nabla^mf)(X_t(x))\Big)\Big|\leq C\|f\|_\infty t^{-\gamma_{k,m}}.\label{ET7}
\end{align}
\et
\begin{proof}
First of all, we assume that
\begin{align}
\mbox{$S_t$ has finite moments of all orders.}\label{Assu}
\end{align}
Under this assumption, it is standard to prove that for any $m,k\in \mN_0$ with $m+k\geq 1$ and $p\geq1$ (cf. \cite[Lemma 3.6]{Zh2}),
\begin{align}
\sup_{x\in\mR^d}\sup_{t\in[0,1]}\mE\left(\|D^m\nabla^k X_t(x)\|_{\mH^{\otimes^m}}^p\right)<+\infty.\label{ET6}
\end{align}
By the chain rule, we have
$$
\nabla^k\mE\Big((\nabla^mf)(X_t(x))\Big)=\sum_{j=1}^k\mE\Big((\nabla^{m+j}f)(X_t(x))G_j(\nabla X_t(x),\cdots,\nabla^k X_t(x))\Big),
$$
where $\{G_j, j=1,\cdots,k\}$ are real polynomial functions.
By Theorem \ref{Th2} and (\ref{ET6}), through cumbersome calculations,
one finds that there exist integer $n$ and $p>1$, $C>0$ such that for all $t\in(0,1)$ and $x\in\mR^d$,
\begin{align*}
\Big|\nabla^k\mE\Big((\nabla^mf)(X_t(x))\Big)\Big|\leq C\|f\|_\infty\mE\left(\|(\det\Sigma_t(x))^{-1}\|^n_p\right)\leq C\|f\|_\infty\|(\det\Sigma_t(x))^{-1}\|^n_{np}.
\end{align*}
Estimate (\ref{ET7}) now follows by Lemma \ref{Le4}.

Now, without assumption (\ref{Assu}), we can use the same argument as in \cite[Section 3.3]{Zh2} to prove (\ref{ET7}).
Since it is completely the same, we omit the details.
\end{proof}

Now, we are in a position to give
\\
\\
{\it Proof of Theorem \ref{Th1}}: 
Notice that for $\nu_S(\dif u)=u^{-1-\frac{\alpha}{2}}\dif u$, {\bf(AS$^\alpha$)} holds.
By estimate (\ref{ET7}) and Sobolev's embedding theorem,
one has (\ref{ET9})  (see \cite[pp.102-103]{Nu}). Equation (\ref{ET8}) follows by It\^o's formula.
The smoothness of $\rho_t(x,y)$ with respect to the time variable $t$ follows
by equation (\ref{ET8}) and a standard bootstrap argument. 

\vspace{5mm}

{\bf Acknowledgements:}

The author is very grateful to Hua Chen, Zhen-Qing Chen, Zhao Dong, Xuhui Peng and Feng-Yu Wang for their quite useful conversations.
This work is supported by NSFs of China (Nos. 11271294, 11325105) and Program for New Century Excellent Talents in University (NCET-10-0654).

%\vspace{5mm}

\end{document}